\documentclass[reqno]{amsart}
\usepackage[a4paper,hmarginratio=1:1]{geometry}

\usepackage{amsmath}
\usepackage{mathabx}
\usepackage{amssymb}
\usepackage{mathtools}
\usepackage{amsthm}
\usepackage{graphicx}
\usepackage{titlesec}
\usepackage{lipsum}
\usepackage{comment}

\usepackage{enumitem}
\usepackage{tikz-cd}
\usepackage{dsfont}
\usepackage{lmodern}
\usepackage{bbm}
\usepackage{adjustbox}
\usepackage{hyperref}
\usepackage{mathrsfs}
\definecolor{dark-red}{rgb}{0.5,0.15,0.15}
\definecolor{dark-blue}{rgb}{0.15,0.15,0.6}
\definecolor{dark-green}{rgb}{0.15,0.6,0.15}
\hypersetup{
    colorlinks, linkcolor=dark-red,
    citecolor=dark-blue, urlcolor=dark-green
}
\usepackage[capitalize]{cleveref}
\newtheorem{theorem}{Theorem}[section]
\newtheorem{lemma}[theorem]{Lemma}
\newtheorem{proposition}[theorem]{Proposition}
\newtheorem{corollary}[theorem]{Corollary}
\theoremstyle{definition}
\newtheorem{definition}[theorem]{Definition}
\newtheorem{example}[theorem]{Example}

\newtheorem{remark}[theorem]{Remark}

\numberwithin{equation}{section}
\newcommand{\stabcat}{\widehat{\textup{Mod}}}

\newcommand{\loc}{\textup{Loc}}

\newcommand{\locm}{\textup{Loc}_{\otimes}}
\newcommand{\colocm}{\textup{Coloc}_{\textup{hom}}}
\newcommand{\coloc}{\textup{Coloc}}

\newcommand{\lhf}{\textup{LH}\mathfrak{F}}

\newcommand{\cat}{\mathfrak{C}}

\titleformat{\section}
  {\normalfont\normalsize\bfseries\filcenter}{\thesection.}{1em}{}
  \titleformat{\subsection}
  {\normalfont\normalsize\bfseries}{\thesubsection.}{1em}{}
\titleformat{\subsubsection}
  {\normalfont\normalsize\bfseries}{\thesubsubsection.}{1em}{}
  \title{(Co)stratification of stable categories for infinite groups}
\author{Gregory Kendall}
\begin{document}
\begin{abstract}
    We classify the localising tensor ideal and colocalising hom-closed subcategories of the stable module category for $\mathrm{LH}\mathfrak{F}$ groups. Along the way we develop techniques to provide similar classifications for other categories of infinite groups.
\end{abstract}\textit{}
\date{\today}
\subjclass[2020]{Primary: 20C07, Secondary: 18G80, 20C12}
\keywords{stable module category, stratification, costratification, infinite groups}
\maketitle
\section{Introduction}
Classifying all modular representations, or even just the finite dimensional representations, over most groups is a wild problem. Instead, it is enlightening to classify representations up to a broader notion of equivalence. Benson, Carlson, and Rickard \cite{bcrthick} show that thick tensor ideal subcategories of the small stable module category of finite dimensional representations for a finite group are in bijection with specialisation closed subsets of the homogeneous prime ideal spectrum of the cohomology ring. This was later extended by Benson, Iyengar, and Krause \cite{repstrat} to the big stable module category of all modules. 
\par 
This fits more widely into the field of tensor triangular geometry \cite{MR2827786}, which unifies such classification theorems across mathematics. Particularly relevant to us is the machinery developed by Benson, Iyengar, and Krause \cite{BIKcol, BIKstrat, BIKloc} in order to classify localising tensor ideals and colocalising hom-closed subcategories in general settings; this culminated in their work on the stable category for finite groups in \cite{repstrat, BIKcol}. 
\par
Barthel, Heard, and Sanders \cite{BHSstrat} have also developed a related notion of stratification for rigidly-compactly generated tensor-triangulated categories. For finite groups, the stable module category is stratified in this sense, see \cite[Theorem~G]{BHSstrat}, which recovers Benson, Iyengar, and Krause's classification of localising tensor ideals. 
\par 
For infinite groups, two main problems arise if we try to apply these techniques to the stable module category, which we constructed in \cite{kendall2024stablemodulecategorymodel}. The first is that the cohomology ring of an infinite group is not necessarily noetherian; this means we cannot apply Benson, Iyengar, and Krause's results since this requires a noetherian graded ring to act on the category. The second issue is that the stable category is not necessarily rigidly-compactly generated. Indeed, there are groups such that the tensor unit is not a compact object, see \cite[Section~5.1]{kendall2024stablemodulecategorymodel}. 
\par 
We can however still classify the subcategories of interest for infinite groups using different methods. In particular, we work with the class of $\textup{LH}\Phi_k$ groups, which contain Kropholler's $\lhf$ groups \cite{krop}. Our main theorem is as follows, where we use $\stabcat(kG)$ to denote the stable category. 
\begin{theorem}\label{introtheorem}(See \Cref{ohyea}) 
     Suppose $k$ is a commutative noetherian ring of finite global dimension and $G$ is an $\textup{LH}\Phi_k$ group. Then there are bijective maps between the following sets.\[\left\{\begin{gathered} \text{Localising tensor ideal}\\ \text{subcategories of $\stabcat(kG)$}
\end{gathered}\;
\right\} \xleftrightarrow{\ \raisebox{-.4ex}[0ex][0ex]{$\scriptstyle{\sim}$}\ }\left\{
\begin{gathered}
  \text{Colocalising hom-closed}\\ \text{subcategories of
    $\stabcat(kG)$}
\end{gathered}
\right\}\; \xleftrightarrow{\ \raisebox{-.4ex}[0ex][0ex]{$\scriptstyle{\sim}$}\ }\left\{
\begin{gathered}
  \text{subsets of}\\ \underset{{F \in \mathcal{A}_{\mathcal{F}}(G)}}{\textup{colim}}\textup{Proj}(H^*(F,k))
\end{gathered}
\right\}\
\]
\end{theorem}
The category $\mathcal{A}_{\mathcal{F}}(G)$ has objects the finite subgroups of $G$ and the morphisms are induced by conjugations by elements of $G$ and subgroup inclusions. For groups where Quillen Stratification is known to hold, for example $\textup{H}_1\mathfrak{F}$ groups of type $\textup{FP}_{\infty}$ \cite{Hennq}, we can further identify the colimit in the above with $\textup{Proj}(H^*(G,k))$. 
\par 
In order to prove \Cref{introtheorem} we first restrict the localising tensor ideals to the finite subgroups, and then apply the classification for finite groups. This method also works with other category, such as $K(\textup{Proj}(kG))$. In fact, we really only need restriction, induction, and coinduction functors which satisfy the Mackey decomposition formula, as well as the projection formula as in \eqref{projmap}, and a suitable version of Chouinard's theorem. 
\par 
We formalise this using Balmer and Dell'Ambrogio's notions of Mackey 2-functors and Green 2-functors \cite{Mackey2,Green2}, with some adaptations to allow our group to be infinite. Working in this context allows us to classify the localising tensor ideals of various different categories using one theory, as well as removing non-canonical choices, such as the choice of representatives for double cosets in the Mackey decomposition formula. 
\subsection*{Acknowledgements}
I am very grateful to my supervisor Peter Symonds for his support and many useful discussions.
\section{Setup}\label{setupsection}
We begin by describing the general setting we require. This is very similar to the setting of Mackey and Green 2-functors \cite{Mackey2,Green2} although we have to make adjustments to include infinite groups. For the necessary 2-categorical background, we refer to \cite[Appendix~A]{Mackey2}. 
\par 
We let $\textbf{Gpd}$ denote the 2-category of (not necessarily finite) groupoids, i.e. categories in which every morphism is invertible. The 1-morphisms are functors and 2-morphisms are given by natural transformations; note that these are necessarily invertible. We will write $\textbf{Gpd}^f$ to denote the sub-2-category where we only consider faithful functors between the groupoids.
\begin{remark}
We can view a group as a groupoid with one object, and in this case a morphism $G_1 \to G_2$ is simply a group homomorphism. Furthermore, a 2-morphism between $f_1,f_2: G_1 \to G_2$ can be identified with an element $g \in G_2$ such that ${}^gf_1 = f_2$, by which we mean that $gf_1(x)g^{-1} = f_2(x)$ for every $x \in G_1$. 
\end{remark}
We then take any 2-full sub-2-category $\mathbb{G}$ of $\textbf{Gpd}$. For example, this could be the sub 2-category of all finite groupoids. We will be particularly interested in the case where $\mathbb{G}$ is the subcategory of all groupoids which can be written as coproducts of certain classes of infinite groups. Note that we don't require $\mathbb{G}$ to be 1-full - for example we may only want to consider the faithful functors between groupoids. 
\par 
We recall the construction of isocomma groupoids as follows. Given two morphisms of groupoids $h:H \to G$ and $k:K \to G$ we have the isocomma groupoid, denoted $(h/k)$, with objects triples $(x,y,g)$ such that $x$ is an object of $H$, $y$ is an object of $K$ and $g: h(x) \to k(y)$ is an invertible morphism in $G$. Morphisms $(x,y,g) \to (x',y'g')$ consist of a pair $(i,u)$ such that $i: x \to x'$ and $u: y \to y'$ are such that $g'h(i) = k(u)g$. This fits into the following 2-cell.
\[\begin{tikzcd}
        & (h/k)\arrow{dl}[swap]{p}\arrow{dr}{q}\\
        H\arrow{dr}[swap]{h} &\overset{\sim}{\Rightarrow}& K\arrow{dl}{k}\\
        & G
    \end{tikzcd}\]
    The 2-cell $\gamma: hp \Rightarrow kq$ is invertible and is given by $(\gamma)_{(x,y,g)} = g$. Note that this 2-cell enjoys a universal property as in \cite[Definition~2.1.1]{Mackey2}.
    \par 
    We also need to recall the notion of mates, see e.g. \cite[Appendix~A.2]{Mackey2} for further information. Let $A,A',B,B'$ be categories and suppose we are given adjunctions $l \dashv r: A \to B$ and $l' \dashv r': A' \to B'$ as well as functors $f_1: B \to B'$ and $f_2:A \to A'$.
    \par 
    There is then a natural bijection $[l'f_2,f_1l] \simeq [f_2r,r'f_1]$, where we use the notation $[-,-]$ to denote the class of natural transformations between the functors. For $\alpha: f_2r \Rightarrow r'f_1$, the image under this bijection is called the left mate of $\alpha$; right mates are defined similarly.
    \par 
    We can describe this bijection explicitly. Let $(\eta, \varepsilon)$ be the unit and counit of the adjunction $l \dashv r$, and similarly for $(\eta',\varepsilon')$. Then the left mate of $\alpha: f_2r \Rightarrow r'f_1$ is the following natural transformation 
    \[\begin{tikzcd}[arrows = Rightarrow]
        l'f_2 \arrow{r}{l'f_2\eta} & l'f_2rl \arrow{r}{l'\alpha l} & l'r'f_1l \arrow{r}{\varepsilon'f_1l} & f_1l
    \end{tikzcd}\]
    For $\beta: l'f_2 \Rightarrow f_1l$ the right mate is defined as follows. 
    \[\begin{tikzcd}[arrows = Rightarrow]
        f_2r \arrow{r}{\eta'f_2r} & r'l'f_2r \arrow{r}{r'\beta r}&r'f_1lr \arrow{r}{r'f_1\varepsilon} & r'f_1
    \end{tikzcd}
    \]
   We can now make our key definition. Let $\textbf{ADD}$ be the 2-category of additive categories and additive functors. We write $\mathbb{G}^{op}$ to denote the 2-category where we have reversed the 1-cells, but not the 2-cells.    
\begin{definition}
A left partial Mackey 2-functor is a strict 2-functor $\cat: \mathbb{G}^{op} \to \textbf{ADD}$ satisfying the following:
\begin{enumerate}[label=(\roman*)]
    \item For any family $\{G_i\}_{i \in I}$ in $\mathbb{G}$, the natural map $\cat(\bigsqcup\limits_{j \in I}G_j) \to \prod\limits_{i \in I}\cat(G_i)$ is an equivalence.
    \item For every faithful functor $i: H \to G$, the restriction $i^*: \cat(G) \to \cat(H)$ preserves coproducts and has a left adjoint $i_!$.
    \item For every iso-comma square
    \[\begin{tikzcd}
        & (i/u)\arrow{dl}[swap]{p}\arrow{dr}{q}\\
        H\arrow{dr}[swap]{i} &\overset{\sim}{\Rightarrow}& K\arrow{dl}{u}\\
        & G
    \end{tikzcd}\]
    there is an isomorphism of the following form
    \[q_!\circ p^* \Rightarrow u^*\circ i_!\]
    which is given by the left mate of $\gamma^*$, where $\gamma: i \circ p \Rightarrow u \circ q$ is the 2-cell in the above iso-comma square. 
\end{enumerate}
We say $\cat$ is a partial Mackey 2-functor if the following two conditions hold, in addition to the above. 
\begin{enumerate}[resume*]
    \item For every faithful functor $i: H \to G$ the restriction $i^*: \cat(G) \to \cat(H)$ has a right adjoint $i_*$
    \item The right mate of $(\gamma^{-1})^*$, where $\gamma$ is as in $(3)$, is an isomorphism $u^* \circ i_* \Rightarrow q_* \circ p^*$
\end{enumerate}
\end{definition}
We will occasionally refer to property $(i)$ as additivity and $(iii),(v)$ as the left and right Beck-Chevalley conditions respectively.
\begin{remark}
    It is worthwhile pointing out the differences between our definition and that of a Mackey 2-functor in \cite{Mackey2}. The key difference is that Balmer and Dell'Ambrogio require that for every faithful $i:H \to G$ there is an isomorphism between the right and left adjoints of the restriction. We are interested in particular in infinite groups, and so we would not expect this in our context. For example, if $H \leq G$ is a subgroup of infinite index, then the usual induction and coinduction of modules will not be isomorphic as functors.
\end{remark}
Consider the 2-category $\textbf{SymMonADD}^{\oplus}$ of symmetric monoidal, idempotent complete additive categories with coproducts, such that the tensor products distribute over the coproducts. In our applications, we will generally be interested in the sub-2-category of tensor-triangulated categories.
\par
\sloppy Suppose further that our partial Mackey 2-functor $\cat$ lifts along the forgetful functor $\textbf{SymMonADD}^{\oplus} \to \textbf{ADD}$ to give a 2-functor $\cat: \mathbb{G}^{op} \to \textbf{SymMonADD}^{\oplus}$. Consider a faithful 1-morphism $h: H \to G$. We obtain the following (left) projection formula for any $X \in \cat(G)$ and $Y \in \cat(H)$.
\begin{equation}\label{projmap}
    h_!(h^*(X) \otimes Y) \to X \otimes h_!(Y)
\end{equation}
There is also a right projection formula, using the right adjoint $h_*$ if it exists, which we will not consider in this work. Therefore, we refer to \eqref{projmap} simply as the projection formula.
\begin{definition}
We say $\cat$ is a left partial Green 2-functor if $\cat$ is a left partial Mackey 2-functor and we have a lifting 
\[\begin{tikzcd}[column sep = 70]
    & \textbf{SymMonADD}^{\oplus}\arrow{d}{\textup{forget}} \\ 
    \mathbb{G}^{op} \arrow[ur,dashed] \arrow{r}{\cat} & \textbf{ADD}
\end{tikzcd}\]
such that the projection maps \eqref{projmap} are all invertible.
\par 
If $\cat$ is a partial Mackey 2-functor and the lifting is to the category of closed symmetric monoidal additive categories with products and coproducts, then we say $\cat$ is a partial Green 2-functor.
\end{definition}
Note that since we are assuming all the tensor products are symmetric, we can exchange the roles of $X$ and $Y$ in the projection formula. 
\par 
We give some examples of (left) partial Green $2$-functors. 
\begin{example}\label{examplemod}
    Any monoidal additive derivator $\mathcal{D}: \textbf{Cat}^{op} \to \textbf{ADD}$ restricts to a left partial Green 2-functor on groupoids cf. \cite[Section~4.1]{Mackey2}, \cite[Theorem~10.1]{Green2}. In particular, $\cat:G \mapsto \textup{Mod}(k)^G$ is a partial Green 2-functor; on groups this takes its value in $\textup{Mod}(kG)$. Similarly, we have a partial Green 2-functor whose value on a group is the category of chain complexes of $kG$-modules.
\end{example}
We can take subfunctors as follows. 
\begin{example}
    Suppose $\cat: \mathbb{G}^{op} \to \textbf{ADD}$ is a partial Green 2-functor. Assume we have for each $G \in \mathbb{G}$ a collection of subcategories $\mathfrak{B}(G) \subseteq \cat(G)$ which are closed under restriction, induction, and coinduction and such that $\mathfrak{B}(G)$ is a tensor ideal in $\cat(G)$. Then $\mathfrak{B}$ inherits the structure of a partial Green 2-functor. 
\end{example}
We can also quotient by such subfunctors. 
\begin{lemma}
    Suppose $\mathfrak{B}$ is a sub-2-functor of $\cat$ as above. Then the additive quotient categories $\cat(G)/\mathfrak{B}(G)$ inherit the structure of a left partial Green 2-functor.
\end{lemma}
\begin{proof}
This follows just as \cite[Proposition~4.2.5]{Mackey2} and \cite[Proposition~10.5]{Green2}. The idea is that by assumption on $\mathfrak{B}(G)$ we know that restriction, induction, and coinduction descend to the quotients as do the adjunctions and Beck-Chevalley conditions. Furthermore, the monoidal structure on the quotients is that inherited from $\cat(G)$, which makes sense since $\mathfrak{B}(G)$ is tensor ideal. The projection formula then holds in the quotient since the adjunctions and monoidal structure descend from $\cat(G)$, where it holds by assumption.
\end{proof}
Note that in the above, we see no reason why the quotient categories should inherit a closed symmetric monoidal structure; if they do then this would obviously give a partial Green 2-functor. 
\par 
Later on, we will focus on situations where the 2-functor takes values in tensor-triangulated categories. We recall some standard terminology for later use. 
\par 
Suppose that $\mathcal{T}$ is a triangulated category with coproducts. An object $X$ is called compact if the natural map $\bigoplus \textup{Hom}_{\mathcal{T}}(X,Y_i) \to \textup{Hom}_{\mathcal{T}}(X,\bigoplus Y_i)$ is an isomorphism for all coproducts $\bigoplus Y_i$. We say that $\mathcal{T}$ is compactly generated by a set $\mathcal{C}$ of compact objects if the smallest localising subcategory containing all $C \in \mathcal{C}$ is the whole of $\mathcal{T}$, where a localising subcategory is a full, triangulated subcategory closed under coproducts. By \cite[Lemma~2.2]{SScomp}, this is equivalent to the statement that if $\textup{Hom}_{\mathcal{T}}(C,X) = 0$ for all $C \in \mathcal{C}$ then $X \cong 0$.
\par 
We say an object $X \in \mathcal{T}$ is rigid (also called strongly dualisable in some sources) if the natural map $\textup{hom}(X,\mathbbm{1})\otimes Y \to \textup{hom}(X,Y)$ is an isomorphism for any $Y \in \mathcal{T}$, where $\textup{hom}(-,-)$ denotes the internal hom in the category. We say that the category is rigidly-compactly generated if it is compactly generated and the compact objects coincide with the rigid objects. We will denote the class of rigid (or equivalently compact) objects by $\cat(H)^d$. The superscript denotes dualisable; we use this notation in order to agree with \cite{barthel2023cosupport}.
\section{Localising, tensor ideals}\label{lcs}
Throughout this section, we fix a 2-full sub-2-subcategory $\mathbb{G}$ of $\textbf{Gpd}$ and a group $G$, which is in the subcategory $\mathbb{G}$. We also assume we have a left partial Green 2-functor $\cat: \mathbb{G}^{op} \to \textbf{SymMonADD}^{\oplus}$. 
\begin{definition}\label{allde}
    Suppose $\mathcal{C}$ is a symmetric monoidal additive category with coproducts. We say that a full subcategory $\mathcal{S}$ of $\mathcal{C}$ is:
\begin{enumerate}[label=(\roman*)]
    \item thick if it is closed under summands
    \item localising if it is thick and closed under coproducts 
    \item tensor ideal if for all $S \in \mathcal{S}$ and $X \in \mathcal{C}$ we have $S \otimes X \in \mathcal{S}$ 
\end{enumerate}
If $\mathcal{C}$ is a tensor-triangulated category, we say $\mathcal{S}$ is a localising tensor ideal if it is a triangulated subcategory, in addition to being localising and tensor ideal as defined above.
\end{definition}
We let $\loc(G)$ be the collection of localising tensor ideals of $\cat(G)$. For any $X \in \cat(G)$ we let $\locm(X)$ be the smallest localising tensor ideal containing $X$. Equivalently, this is the intersection of all localising tensor ideals which contain $X$. 
\par 
Our aim is to reduce the question of classifying the localising, tensor ideal subcategories of $\cat(G)$ to a collection of subgroups of $G$. Instead of a collection of subgroups of $G$, it is more convenient in our context to take a subcategory of the comma 2-category $\mathbf{Gpd}^f_{/G}$, the construction of which we now recall as in \cite[Definition~A.1.21]{Mackey2}. 
\begin{definition}
Let $G$ be a groupoid, then the comma 2-category $\mathbb{G}_{/G}$ has objects $(H,i_H)$ where $H$ is a groupoid and $i_H: H \to G$ is a 1-morphism in $\mathbb{G}$. A 1-morphism $(H,i_H) \to (K,i_K)$ is a pair $(f,\gamma_f)$ such that $f: H \to K$ is a 1-morphism in $\mathbb{G}$ and $\gamma_f: i_Kf \Rightarrow i_H$ is a 2-morphism in $\mathbb{G}$. Finally, a 2-morphism $(f,\gamma_f) \Rightarrow (f', \gamma_{f'})$ is a 2-morphism $\alpha: f \to f'$ of $\mathbb{G}$ such that $\gamma_g(i_K\alpha) = \gamma_f$.
\end{definition} 
We now take a 1-full, 2-full 2-subcategory $\mathfrak{H}$ of $\mathbf{Gpd}^f_{/G}$, which we think of as a generalised collection of subgroups of $G$. Furthermore, we assume that $\mathfrak{H}$ is closed under iso-comma squares in the following sense. Suppose  $(H,h),(K,k) \in \mathfrak{H}$ and form the iso-comma square 
 \[\begin{tikzcd}
        & (h/k)\arrow{dl}[swap]{p}\arrow{dr}{q}\\
        H\arrow{dr}[swap]{h} &\overset{\sim}{\Rightarrow}& K\arrow{dl}{k}\\
        & G
    \end{tikzcd}\]
We require that $((h,k),hp) \in \mathfrak{H}$. Since $\mathfrak{H}$ is 2-replete, we can equivalently require $((h,k),kq) \in \mathfrak{H}$.
\par 
 We start by proving a crucial lemma.
 \begin{lemma}\label{summandlem}
     Suppose we have a faithful $1$-morphism $h: H \to G$. For any $X \in \cat(H)$, we have that $X$ is a summand of $h^*h_!(X)$. 
 \end{lemma}
 \begin{proof}
 It is shown in \cite[Proposition~3.1.3]{Mackey2} that there exists a fully faithful embedding $\Delta: H \to (h/h)$. As in the proof of \cite[Proposition~3.2.1]{Mackey2} this induces an equivalence 

\[\begin{tikzcd}[column sep = 80, row sep = 30]\cat((h/h)) \arrow{r}{(\Delta^* \hspace{0.5cm}
 (\textup{inc}_C^*)_C)} &\cat(H) \oplus \prod_{C \cap \Delta(H) = \emptyset}\cat(C)\end{tikzcd}\]
 where $C$ runs over the connected components of $(h/h)$ disjoint from $\Delta(H)$. 
 \par 
 Take any $Y, Z \in \cat((h,h))$ and consider the following natural isomorphisms. 
 \begin{align}
     \cat((h/h))(Y,Z) & \cong \cat(H)(\Delta^*(Y),\Delta^*(Z)) \oplus \prod\limits_{C}\cat(C)(\textup{inc}_C^*(Y),\textup{inc}_C^*(Z)) \\ & \cong \cat((h,h))(\Delta_!\Delta^*(Y),Z) \oplus \prod\limits_{C}\cat((h/h))(\textup{inc}_{C!}\textup{inc}_C^*(Y),Z) \\ &\cong \cat((h/h))(\Delta_!\Delta^*(Y) \oplus \bigoplus\limits_{C}\textup{inc}_{C!}\textup{inc}_C^*(Y),Z)
 \end{align}
 The Yoneda Lemma now implies that $\Delta_!\Delta^*(Y)$ is a summand of $Y$.
 \par 
    We let $p,q: (h/h) \to H$ be the two natural projections. From the Beck-Chevalley condition, we know that the mate $q_!p^* \Rightarrow h^*h_!$ is an isomorphism. Therefore, it suffices to show that $X$ is a summand of $q_!p^*(X)$ for any $X \in \cat(H)$. We have shown that $q_!\Delta_!\Delta^*p^*(X)$ is a summand of $q_!p^*(X)$, and so we only need to show that $q_!\Delta_!\Delta^*p^*(X) \cong X$. However, this follows since $p\Delta = q\Delta = \textup{Id}_H$, see \cite[Proposition~3.1.3]{Mackey2}.
 \end{proof}
 We will also use the following repeatedly; it is a direct consequence of the projection formula. 
 \begin{lemma}\label{profor}
     Suppose we have a faithful $1$-morphism $h:H \to G$. Then for any $X \in \cat(G)$ we have $h_!h^*(X) \cong X \otimes h_!(\mathbbm{1}_H)$.
 \end{lemma}
 \begin{proof}
     Since $h^*$ is a symmetric monoidal functor, it preserves the tensor unit. Therefore, we have the following isomorphisms. 
     \[h_!h^*(X) \cong h_!(h^*(X) \otimes h^*(\mathbbm{1}_G)) \cong h_!(h^*(X) \otimes \mathbbm{1}_H) \cong X \otimes h_!(\mathbbm{1}_H)\]
     The final isomorphism is simply the projection formula.
 \end{proof}
Suppose we have a faithful 1-morphism $h:H \to G$ and a subcategory $\mathcal{S}$ of $\cat(H)$. We write $(h^*)^{-1}(\mathcal{S})$ to be the full subcategory of $\cat(G)$ of all objects $X$ such that $h^*(X) \in \mathcal{S}$; we also define $h_!^{-1}$ in the obvious way. The following says that these functions preserve the properties of being localising tensor ideals. 
\begin{proposition}\label{bothdefined}
    Let $h: H \to G$ be a faithful 1-morphism in $\textbf{Gpd}$. Then $(h^{*})^{-1}:\loc(H) \to \loc(G)$ is well defined. Similarly, $h_!^{-1}$ is a well defined function $\loc(G) \to \loc(H)$. 
\end{proposition}
\begin{proof}
    Given that $h^*$ is symmetric monoidal and preserves coproducts, by assumption, this part is straightforward. 
    \par
    We prove the claim about the left adjoint $h_!$. Being a left adjoint, it preserves coproducts and so clearly if $\mathcal{L} \in \loc(G)$ then $h_!^{-1}(\mathcal{L})$ is closed under coproducts and summands. 
    \par 
    We need to show that $h_!^{-1}(\mathcal{L})$ is closed under tensoring, and so take $X \in h_!^{-1}(\mathcal{L})$ and $Y \in \cat(H)$. We have the following isomorphism from the projection formula. 
    \[h_!(h^*h_!(X) \otimes Y) \cong h_!(X) \otimes h_!(Y)\]
    It follows that $h^*h_!(X) \otimes Y \in h_!^{-1}(\mathcal{L})$. We have assumed that the tensor products commutes with coproducts, and so by \Cref{summandlem} we know that $X \otimes Y$ is a summand of $h^*h_!(X) \otimes Y$; this completes the proof since we have already shown that $h_!^{-1}(\mathcal{L})$ is closed under summands. 
\end{proof}
We have another description of the function $h_!^{-1}$, which is perhaps more natural to consider. In the following, we use the notation $h^*(\mathcal{L})$ to denote the collection of all objects of the form $h^*(X)$ for $X \in \mathcal{L}$.
\begin{lemma}\label{lreqca}
    Let $h: H \to G$ be a faithful 1-morphism in $\textbf{Gpd}$ and take $\mathcal{L} \in \loc(G)$. Then $h_!^{-1}(\mathcal{L}) = \locm(h^*(\mathcal{L}))$.
\end{lemma}
\begin{proof}
    Take any $X \in \mathcal{L}$. The projection formula shows that $h^*(X) \in h_!^{-1}(\mathcal{L})$. We know from \Cref{bothdefined} that $h_!^{-1}(\mathcal{L})$ is a localising tensor ideal, and so we conclude that $\locm(h^*(\mathcal{L})) \subseteq h_!^{-1}(\mathcal{L})$. 
    \par 
    On the other hand, suppose $Y \in \cat(H)$ is such that $h_!(Y) \in \mathcal{L}$. Consider the localising tensor ideal $\mathcal{J} = (h^*)^{-1}(\locm(h^*(\mathcal{L})))$. Clearly, $\mathcal{L} \subseteq \mathcal{J}$ and hence we see that $h^*h_!(Y) \in \locm(h^*(\mathcal{L}))$. By \Cref{summandlem} we know that $Y$ is a summand of $h^*h_!(Y)$ and hence also we have $Y \in \locm(h^*(\mathcal{L}))$. This shows the desired equality.
\end{proof}
Note that there is a priori no reason for there to be a set (i.e. not a proper class) of localising tensor ideals in $\cat(G)$. However, from now on we make the assumption that $\loc(H)$ is a set for each $(H,h) \in \mathfrak{H}$; this allows us to make the following key definition.
\begin{definition}\label{ldefinition}
We now define a 2-functor $\mathbb{L}: \mathfrak{H} \to \mathbf{Set}$; this will depend on $\cat$, although we don't include that in the notation in order to avoid overcrowding. 
\par 
On objects we take a pair $(H,h)$ and let $\mathbb{L}((H,h))$ be $\loc(\cat(H))$. 
\par 
Given a 1-morphism $(H,h) \to (K,k)$, i.e. a pair $(f, \gamma_f)$, we let $\mathbb{L}((f,\gamma_f)) = f_!^{-1}$. This is well-defined by \Cref{bothdefined}.
\par 
For 2-morphisms we have no choice since the only invertible 2-morphisms in $\textbf{Set}$ are equalities. We must check that given $(f,\gamma_f) \Rightarrow (f',\gamma_{f'})$, i.e. an $\alpha: f \Rightarrow f'$ such that the appropriate diagrams commute, we have $(f_!)^{-1} = (f'_!)^{-1}$. We know that $\cat(\alpha): f^* \Rightarrow f'^*$ is a natural isomorphism and hence $f_!$ and $f'_!$ are naturally isomorphic. Therefore, since localising tensor ideals are closed under isomorphisms, it follows that there is an equality $(f_!)^{-1} = (f'_!)^{-1}$ as desired.
\end{definition}
It is a straightforward verification that $\mathbb{L}$ defines a strict 2-functor. 
\par 
We can calculate the 2-limit of $\mathbb{L}$ as a limit in the 1-category $\textup{Set}$. This is because if an object satisfies the universal property for a 2-limit, it must also satisfy the universal property for the 1-limit. Explicitly, we have 
\begin{equation}\label{limite}
    \lim\limits_{(H,h) \in \mathfrak{H}}\mathbb{L}((H,h)) \cong \{ (\mathcal{L}_H)_{(H,h) \in \mathfrak{H}} \in \prod\limits_{(H,h) \in \mathfrak{H}}\loc(\cat(H))\ |\ \textup{for all }(f,\gamma_f):(H,h) \to (K,k), \mathcal{L}_{H} = f_!^{-1}(\mathcal{L}_{K})\}
\end{equation}
For any object $(H,h)$ of $\mathfrak{H}$ we have a map $(h_!)^{-1}: \loc(G) \to \loc(H)$. Suppose we have a 1-morphism $(f,\gamma_f): (H,h) \to (K,k)$ in $\mathfrak{H}$. Since $\gamma_f$ is a 2-isomorphism, we know that $\cat(\gamma_f): f^*k^{*} \Rightarrow h^*$ is a natural isomorphism and consequently we have another natural isomorphism $k_!f_! \cong h_!$. This shows that as maps of sets we have an equality $f_!^{-1}k^{-1}_! = h_!^{-1}$. Altogether, these combine to give a natural map $\loc(G) \to \lim\limits_{(H,h) \in \mathfrak{H}}\mathbb{L}((H,h))$.
\par
We are interested in determining when this natural map is an isomorphism. In particular, this would imply that $\loc(G)$ is in fact a set. Under our assumptions so far, we can show that it is always surjective by exhibiting a right inverse. 
\par 
For this we introduce the following notation. 
\begin{definition}
    Consider an element $(\mathcal{L}_H) \in \lim\limits_{(H,h) \in \mathfrak{H}}\mathbb{L}((H,h))$. We define the map 
\[\bigwedge\limits_{(H,h) \in \mathfrak{H}}(h^{*})^{-1}: (\mathcal{L}_H) \mapsto \bigcap\limits_{(H,h) \in \mathfrak{H}} (h^*)^{-1}(\mathcal{L}_H)\]
\end{definition}
This means that $X \in \bigwedge\limits_{(H,h) \in \mathfrak{H}}(h^{*})^{-1}(\mathcal{L}_H)$ if and only if $h^*(X) \in \mathcal{L}_H$ for each $(H,h) \in \mathfrak{H}$. We use this notation since we are essentially taking the meet of all the localising ideals $(h^*)^{-1}(\mathcal{L}_H)$ in the lattice $\loc(G)$; see \Cref{framsection}.
\begin{proposition}\label{rightinverse}
    Suppose $\mathfrak{H}$ is a 1-full, 2-full 2-subcategory of $\textbf{Gpd}^f_{/G}$ closed under iso-comma squares. 
    The map $\bigwedge\limits_{(H,h) \in \mathfrak{H}}(h^{*})^{-1}$ is a right inverse to the natural map $\loc(G) \to\lim\limits_{(H,h) \in \mathfrak{H}}\mathbb{L}((H,h))$. 
    \par 
    Explicitly, let $(\mathcal{L}_H)$ be an element of $\lim\limits_{(H,h) \in \mathfrak{H}}\mathbb{L}((H,h))$ as in \eqref{limite}. Then for any $(K,k) \in \mathfrak{H}$ we have $k_!^{-1}\bigwedge\limits_{(H,h) \in \mathfrak{H}} (h^{*})^{-1}(\mathcal{L}_H) = \mathcal{L}_K$. 
\end{proposition}
\begin{proof}
    Let $(\mathcal{L}_H)$ be an element of the limit. We need to show that for any $(K,k) \in \mathfrak{H}$ we have $k_!^{-1}\bigwedge\limits_{(H,h) \in \mathfrak{H}} (h^{*})^{-1}(\mathcal{L}_H) = \mathcal{L}_K$. 
    \par
    Suppose that $X \in k_!^{-1}\bigwedge\limits_{(H,h) \in \mathfrak{H}} (h^{*})^{-1}(\mathcal{L}_H)$. In particular, $k^*k_!(X) \in \mathcal{L}_K$. From \Cref{summandlem} we know that $\textup{id}_K$ is a summand of $k^*k_!$ and hence $X \in \mathcal{L}_K$; altogether this shows the inclusion $k_!^{-1} \bigwedge\limits_{(H,h) \in \mathcal{A}(G)} (h^{*})^{-1}(\mathcal(\mathcal{L}_H)) \subseteq \mathcal{L}_K$.
    \par
    For the reverse inclusion, we show that if $X \in \mathcal{L}_K$ then $h^*k_!X \in \mathcal{L}_H$ for any $H \in \mathfrak{H}$. So, fix some $H$ and consider the following iso-comma square. 
    \[\begin{tikzcd}
        & (k/h)\arrow{dl}[swap]{p}\arrow{dr}{q}\\
        K\arrow{dr}[swap]{k} &\overset{\sim}{\Rightarrow}& H\arrow{dl}{h}\\
        & G
    \end{tikzcd}\] 
    By the Beck-Chevalley condition, we know that the mate $q_!p^* \Rightarrow h^*k_!$ is an isomorphism. 
    Therefore, it suffices to show that $q_!p^*(X) \in \mathcal{L}_H$. By the projection formula, see \Cref{profor}, we have an isomorphism $p_!p^*(X) \cong X \otimes p_!(\mathbbm{1}_{(k/h)})$.
    \par 
    As $\mathcal{L}_K$ is tensor ideal, it now follows that $X \otimes p_!(\mathbbm{1}_{(k/h)}) \cong p_!p^*X \in \mathcal{L}_K$ and so $p^*X \in p_!^{-1}\mathcal{L}_K$. 
    \par 
    Since we started with an element of the limit, we know that $p_!^{-1}\mathcal{L}_K = \mathcal{L}_{(k/h)} = q_!^{-1}\mathcal{L}_H$; here we use the assumption that $\mathfrak{H}$ is closed under isocomma squares. It follows that $p^*X \in  q_!^{-1}\mathcal{L}_H$, i.e. $q_!p^*X \in \mathcal{L}_H$ as desired.
    \end{proof}
The above works quite generally; however to exhibit a left inverse we need to impose some conditions on our group $G$. 
\begin{definition}\label{detecteddefinition}
    We say that localising tensor ideals for $G$ are detected on restriction to $\mathfrak{H}$ if for every $X \in \cat(G)$ we have 
    \[\locm(X) = \locm\{\bigoplus\limits_{(H,h) \in \mathfrak{H}}h_!h^*(X)\}\]
\end{definition}
The following is then the main result of this section. 
\begin{theorem}\label{bigboy}
    Suppose we have a 1-full, 2-full 2-subcategory $\mathfrak{H}$ of $\mathbf{Gpd}^f_{/G}$ which is closed under isocomma squares and is such that localising ideals are detected on restriction to $\mathfrak{H}$. Then the natural map $\loc(G) \to\lim\limits_{(H,h) \in \mathfrak{H}}\mathbb{L}((H,h))$ is invertible.
\end{theorem}
\begin{proof}
    From \Cref{rightinverse}, we know that the natural map has a right inverse given by $\bigwedge\limits_{(H,h) \in \mathfrak{H}}(h^{*})^{-1}$. We must show this is also a left inverse. For this, take some $\mathcal{L} \in \loc(G)$ and we need to show that $\bigwedge\limits_{(H,h) \in \mathfrak{H}}(h^{*})^{-1}h_!^{-1}(\mathcal{L}) = \mathcal{L}$. Equivalently, we need to show that for any $X \in \cat(G)$ we have $X \in \mathcal{L} \iff h_!h^*(X) \in \mathcal{L}$ for all $(H,h) \in \mathfrak{H}$. However, this follows by assumption since localising ideals are detected on $\mathfrak{H}$.
\end{proof}
\subsection{Reducing to subgroups}
It is obviously important to understand situations in which we have a subcategory $\mathfrak{H}$ such that localising tensor ideals for $G$ are detected on restriction to $\mathfrak{H}$. In general, we will have a collection of subgroups of $G$ and we want to reduce the classification of localising ideals to such subgroups. However, in order to apply \Cref{bigboy} we require a subcategory $\mathfrak{H}$ as in the statement, and not just a collection of subgroups. We therefore make the following definition.
\begin{definition}\label{frakde}
 Suppose we have a collection $\mathcal{H}$ of subgroups of $G$ which is closed under conjugations and subgroups. We define $\mathfrak{H}_{\mathcal{H}}$ to be the 1-full, 2-full subcategory of $\textbf{Gpd}^f_{/G}$ of objects $(K,k)$ such that $k$ factors through the inclusion $H \to G$ for some $H \in \mathcal{H}$\end{definition}
 We may view $\mathcal{H}$ as a subcategory of $\mathfrak{H}_{\mathcal{H}}$, since each $H \in \mathcal{H}$ comes with an inclusion $h: H \to G$ and then clearly $(H,h) \in \mathfrak{H}_{\mathcal{H}}$. We will abuse notation and continue to refer to this subcategory as $\mathcal{H}$.
\begin{lemma}\label{isocommaclosure}
    The subcategory $\mathfrak{H}_{\mathcal{H}}$ is closed under isocomma squares.
\end{lemma}
\begin{proof}
    Suppose $(H,h), (K,k) \in \mathfrak{H}_{\mathcal{H}}$ and consider the isocomma square as follows. 
    \[\begin{tikzcd}
        & (h/k)\arrow{dl}[swap]{p}\arrow{dr}{q}\\
        H\arrow{dr}[swap]{h} &\Rightarrow& K\arrow{dl}{k}\\
        & G
    \end{tikzcd}\]
   By definition, there exists some $F \in \mathcal{H}$ and inclusion $f: F \to G$ such that $h$ factors through $f$. But then $hp$ also clearly factors through $F$, showing that $((h/k),hp) \in \mathfrak{H}_{\mathcal{H}}$. 
\end{proof}
\begin{proposition}\label{NCE}
    Suppose $\mathcal{H}$ is a collection of subgroups of $G$, and let $h:H \to G$ denote the inclusion for each $H \in \mathcal{H}$. Assume furthermore that for any $X \in G$ we have \[\locm(X) = \locm\{\bigoplus\limits_{H \in \mathcal{H}}h_!h^*(X)\}\]  Then localising ideals are detected on restriction to $\mathfrak{H}_{\mathcal{H}}$ and so the natural map $\loc(G) \to \lim\limits_{(K,k) \in \mathfrak{H}_{\mathcal{H}}}\mathbb{L}((K,k))$ is invertible.
\end{proposition}
\begin{proof}
    We will view $\mathcal{H}$ along with the inclusions $h:H \to G$ for $H \in \mathcal{H}$ as a subcategory of $\mathfrak{H}_{\mathcal{H}}$. It follows immediately that $\locm\{\bigoplus\limits_{H \in \mathcal{H}}h_!h^*(X)\} \subseteq \locm\{\bigoplus\limits_{(K,k) \in \mathfrak{H}_{\mathcal{H}}}k_!k^*(X)\}$.
    \par 
    On the other hand, the projection formula in the form of \Cref{profor} shows that $k_!k^*(X) \in \locm(X)$ for every $(K,k) \in \mathfrak{H}_{\mathcal{H}}$. Hence, $\locm\{\bigoplus\limits_{(K,k) \in \mathfrak{H}_{\mathcal{H}}}k_!k^*(X)\} \subseteq \locm(X)$. However, by assumption we know that $\locm(X) = \locm\{\bigoplus\limits_{H \in \mathcal{H}}h_!h^*(X)\}$, completing the proof that localising ideals are detected on restriction to $\mathfrak{H}_{\mathcal{H}}$. 
    \par 
    The statement about the invertibility of the map now follows from \Cref{bigboy}, which we can apply thanks to \Cref{isocommaclosure}.
\end{proof}
In fact, we can rewrite this limit over a more familiar category. Before this, we recall that the underlying category of a (2,1)-category $\textbf{C}$, denoted $\textbf{C}^{(1)}$, is the category where objects are those in $\textbf{C}$ and morphisms are the $1$-morphisms of $\textbf{C}$. Since $\textbf{C}$ is a strict 2-category, it makes sense to take the compositions of $\textbf{C}^{(1)}$ to just be those induced by $\textbf{C}$.
\par
Define a 1-category $\mathcal{A}_{\mathcal{H}}(G)$ as follows. The objects are the subgroups $H \in \mathcal{H}$, and morphisms are induced by conjugation by elements of $G$ and inclusions. Then $\mathfrak{H}_{\mathcal{H}}^{(1)}$ is essentially an extension of $\mathcal{A}_{\mathcal{H}}(G)$ where we have added in coproducts. We make this precise in the following. 
\begin{lemma}\label{qle}
Consider an object $(H,h)$ of $\mathfrak{H}_{\mathcal{H}}$. Then $H$ is isomorphic, as a groupoid, to $\bigsqcup\limits_{i \in I}H_i$ such that each $H_i \in \mathcal{H}$. 
\end{lemma}
\begin{proof}
    We can write $H$ as a coproduct of one object groupoids, say $\bigsqcup\limits_{j \in J}K_j$. By definition, the faithful functor $h: H \to G$ factors through the inclusion $f:F \to G$ for some $F \in \mathcal{H}$ and so we have a faithful functor $H \to F$. For each $K_j$, this gives a (necessarily faithful) functor $K_j \to F$. Therefore, we can view $K_j$ as a subgroup of $F$. Since $\mathcal{H}$ is assumed to be closed under subgroups, the result then follows.
\end{proof}
\begin{lemma}
    Suppose we have a collection $\{\mathcal{T}_i\}_{i \in I}$ of additive categories with coproducts. The coproduct in the product $\prod\limits_{i \in I}\mathcal{T}_i$ is computed component wise, i.e. given a collection of objects $(X_i^j) \in \prod\limits_{i \in I}\mathcal{T}_i$ for $j \in J$, then 
    \[\textup{Hom}_{\prod\mathcal{T}_i}((\bigoplus\limits_{j \in J}X_i^j)_i, (Z_i)_i) \cong \prod\limits_{j \in J}\textup{Hom}_{\prod\mathcal{T}_i}((X_i^j)_{i},(Z_i)_i)\] for any $(Z_i) \in \prod\mathcal{T}_i$. 
\end{lemma}
\begin{proof}
    This follows from the following isomorphisms
    \begin{align}
        \textup{Hom}_{\prod\mathcal{T}_i}((\bigoplus\limits_{j \in J}X_i^j)_{i \in I}, (Z_i)) &= \prod\limits_{i \in I}\textup{Hom}_{\mathcal{T}_i}(\bigoplus\limits_{j \in J}X_i^j,Z_i)\\ &\cong \prod\limits_{i \in I}\prod\limits_{j \in J}\textup{Hom}_{\mathcal{T}_i}(X_i^j,Z_i) \\ &\cong \prod\limits_{j \in J}\prod\limits_{i \in I}\textup{Hom}_{\mathcal{T}_i}(X_i^j,Z_i)\\ &= \prod\limits_{j \in J}\textup{Hom}_{\prod\limits_{i \in I}\mathcal{T}_i}((X_i^j)_i,(Z_i)_i)
    \end{align}
    The first and final equalities are by definition of the morphisms in the product of categories. The isomorphisms are by the definition of the coproduct in each $\mathcal{T}_i$, and the fact that products commute with products. 
\end{proof}
In the following, we let $\loc(\mathcal{T})$ denote the set of localising tensor ideals of $\mathcal{T}$. 
\begin{lemma}\label{propres}
    Suppose we have a product $\prod\limits_{i \in I} \mathcal{T}_i$ of symmetric monoidal categories with coproducts. Then $\loc(\prod\limits_{i \in I}\mathcal{T}_i) \cong \prod\limits_{i \in I}\loc(\mathcal{T}_i)$.
\end{lemma}
\begin{proof}
    For each $j \in I$ let $p_j:\prod\limits_{i \in I}\mathcal{T}_i \to \mathcal{T}_j$ be the projection functor. Take any $\mathcal{L} \in \loc(\prod\limits_{i \in I}\mathcal{T}_i)$ and consider $p_j(\mathcal{L})$. First note that this is a localising tensor ideal of $\mathcal{T}_j$: this is easy to see since everything is computed pointwise in the product.
    \par 
    These projection maps assemble to give a map $\loc(\prod\limits_{i \in I}\mathcal{T}_i) \to \prod\limits_{i \in I}\loc(\mathcal{T}_i)$, i.e. $\mathcal{L} \mapsto \prod\limits_{i \in I}p_i(\mathcal{L})$. We claim this map is an isomorphism. 
    \par 
    We begin by showing injectivity, and so suppose $\mathcal{L},\mathcal{J} \in \loc(\prod\limits_{i \in I}\mathcal{T}_i)$ are such that $p_i(\mathcal{L}) = p_i(\mathcal{J})$ for each $i \in I$. We must show that $\mathcal{L} = \mathcal{J}$, and so take some object $(X_i)_{i \in I}\in \mathcal{L}$. 
    \par 
    Take any $j \in I$. Then $X_j = p_j((X_i)_{i \in I}) \in p_j(\mathcal{L})$ by definition. Since we have assumed $p_j(\mathcal{L}) = p_j(\mathcal{J})$ we see that $X_j \in p_j(\mathcal{J})$. This means, again by definition, that there exists some object $(Y_i)_{i \in I} \in \mathcal{J}$ such that $Y_j = X_j$. 
    \par 
    We consider the object $\mathcal{X}_j := (\dots,0, X_j, 0, \dots) \in \prod\limits_{i \in I}\mathcal{T}_i$ which has the zero object everywhere except in index $j$ where it has $X_j = Y_j$. Since coproducts are taken pointwise, we see that this object is a summand of $(Y_i)_{i \in I}$. Therefore, $\mathcal{X}_j \in \mathcal{J}$.
    \par 
    By definition, $\bigoplus\limits_{i \in I}\mathcal{X}_i = (X_i)_{i \in I}$. We know $\mathcal{J}$ is closed under coproducts, and hence $(X_i)_{i \in I} \in \mathcal{J}$. This shows that $\mathcal{L} \subseteq \mathcal{J}$; the reverse inclusion follows by the same argument, with the roles of $\mathcal{L}$ and $\mathcal{J}$ reversed. 
    \par 
    We now turn to surjectivity. Take an object $\prod\limits_{i \in I}\mathcal{L}_i \in \prod\limits_{i \in I}\textup{Loc}(\mathcal{T}_i)$.  Consider the full subcategory $\mathcal{L}$ of $\prod\limits_{i \in I}\mathcal{T}_i$ of elements $(X_i)_{i \in I}$ such that $X_i \in \mathcal{L}_i$ for each $i \in I$. This is easily checked to be a localising tensor ideal and, by construction, for any $j \in I$ we have $p_j(\mathcal{L}) = \mathcal{L}_i$. 
\end{proof}
For the following, we abuse notation slightly: our 2-functor $\mathbb{L}$ can be viewed as a 1-functor $\mathcal{A}_{\mathcal{H}}(G) \to \textup{Set}$, which we continue to call $\mathbb{L}$. We may also view $\mathbb{L}$ as a functor $\mathfrak{H}_{\mathcal{H}}^{(1)} \to \textup{Set}$ by precomposing with the natural inclusion functor $\mathfrak{H}_{\mathcal{H}}^{(1)} \to \mathfrak{H}_{\mathcal{H}}$.
\begin{proposition}\label{limitiso}
    There is an isomorphism $\lim\limits_{H \in \mathcal{A}_{\mathcal{H}}(G)}\mathbb{L}(H) \cong \lim\limits_{(H,h) \in \mathfrak{H}_{\mathcal{H}}^{(1)}}\mathbb{L}((H,h))$
\end{proposition}
\begin{proof}
    By \Cref{qle} we see that every object in $\mathfrak{H}_{\mathcal{H}}$ is equivalent to a coproduct of subgroups in $\mathcal{H}$ along with their inclusions into $G$. Consider such a coproduct $\bigsqcup\limits_{i \in I}H_i$. Then \[\mathbb{L}(\bigsqcup\limits_{i \in I}H_i) = \textup{Loc}(\cat(\bigsqcup\limits_{i \in I}H_i)) = \textup{Loc}(\prod\limits_{i \in I}\cat(H_i)) = \prod\limits_{i \in I}\textup{Loc}(\cat(H_i)) = \prod\limits_{i \in I}\mathbb{L}(H_i)\]
    The first and final equalities are by definition. The second equality holds because $\cat$ is a partial Green 2-functor, and the third equality follows from \Cref{propres}.
    \par 
    This means that $\mathbb{L}$ sends coproducts to products. We know from \Cref{qle} that $\mathcal{A}_{\mathcal{H}}(G)$ and $\mathfrak{H}_{\mathcal{H}}^{(1)}$ differ only by coproducts, and so we can see that the limit over these two categories will be the same.
\end{proof}
\begin{corollary}\label{needtobeframe}
    Suppose we have a collection of subgroups $\mathcal{H}$ of $G$, with inclusions $h:H \to G$ for each $H \in \mathcal{H}$. Furthermore, suppose that \[\locm(X) = \locm\{\bigoplus\limits_{H \in \mathcal{H}}h_!h^*(X)\}\]
    Then the natural map $\loc(G) \to \lim\limits_{H \in \mathcal{A}_{\mathcal{H}}(G)}\mathbb{L}(H)$ is an isomorphism.
\end{corollary}
\begin{proof}
    From \Cref{NCE} we know that the natural map $\loc(G) \to \lim\limits_{(H,h) \in \mathfrak{H}_{\mathcal{H}}}\mathbb{L}((H,h))$ is invertible. Since we are taking the limit to the discrete 2-category $\textup{Set}$, we may identify $\lim\limits_{(H,h) \in \mathfrak{H}_{\mathcal{H}}}\mathbb{L}((H,h))$ with $\lim\limits_{(H,h) \in \mathfrak{H}_{\mathcal{H}}^{(1)}}\mathbb{L}((H,h))$ cf. \cite[Proposition~1.23]{MAILLARD2022108187}. This is in turn isomorphic to $\lim\limits_{H \in \mathcal{A}_{\mathcal{H}}(G)}\mathbb{L}(H)$ by \Cref{limitiso}.
\end{proof}
\subsection{The detection property}\label{subsecd}
We now turn to conditions which allow us to have a collection of subgroups $\mathcal{H}$ as in \Cref{needtobeframe}. For the rest of this section, we assume that $\cat$ takes values in tensor-triangulated categories.
\begin{lemma}\label{Alem}
Assume $\cat(G)$ is a tensor-triangulated category with tensor unit $\mathbbm{1}$. 
    Suppose $A$ is an object of $\cat(G)$ such that $\locm(A) = \cat(G)$. Then for any $\mathcal{L} \in \loc(G)$ and $X \in \cat(G)$ we have $X \in \mathcal{L} \iff X \otimes A \in \mathcal{L}$.
\end{lemma}
\begin{proof}
    The forward implication follows immediately since $\mathcal{L}$ is closed under tensoring with arbitrary objects by definition. On the other hand, suppose $X \otimes A \in\mathcal{L}$ and consider the full subcategory of $\cat(G)$ consisting of objects $Y$ such that $X \otimes Y \in \mathcal{L}$. It is easy to see this is a  localising tensor ideal. By assumption, it contains $A$ and hence contains $\locm(A) = \cat(G)$; from this it follows that $X \cong X \otimes \mathbbm{1} \in\mathcal{L}$. 
\end{proof} 
\begin{definition}
    Suppose we have a collection $\mathcal{H}$ of subgroups of $G$, which is closed under taking intersection and conjugation. We say that the detection property holds if for any object $X \in \cat(G)$: 
\begin{equation}\label{dp}
    X \cong 0 \iff h^*(X) \cong 0 \textup{ for all } H \in \mathcal{H}
\end{equation}
where $h: H \to G$ is the inclusion.
\end{definition}
\begin{example}
    Let $G$ be a finite group, $k$ a field of characteristic $p > 0$ and consider the Mackey 2-functor which takes values in the stable module category. The detection property then holds for $\mathcal{H}$ the collection of elementary abelian $p$-subgroups by Chouinard's theorem. 
\end{example}
For the following, recall the definition of compact objects and compactly generated triangulated categories from \Cref{setupsection}.
\begin{lemma}\label{cop2}
    Suppose $\cat(H)$ is a compactly generated triangulated category for each $H \in \mathcal{H}$, with a set of compact generators $\mathcal{G}_H$, and assume that the detection property holds for $\mathcal{H}$. Then $\mathcal{G} = \{h_!X\ |\ X \in \mathcal{G}_H \textup{ for some } H \in \mathcal{H}\}$ is a set of compact generators for $\cat(G)$.
\end{lemma}
\begin{proof}
    Suppose $Y \in \cat(G)$ is such that $\textup{Hom}_{\cat(G)}(\mathcal{G},Y) = 0$. Then for any $H$ and any $X \in \mathcal{G}_H$ we have, by adjunction, that $\textup{Hom}_{\cat(H)}(X, h^*Y) = 0$. Since $\mathcal{G}_H$ are a set of compact generators, it follows that $h^*Y \cong 0$ for every $H$. Using the detection property we conclude that $Y \cong 0$. This implies that $\mathcal{G}$ is a set of compact generators by \cite[Lemma~2.2]{SScomp}.
\end{proof}
\begin{theorem}\label{projloclem}
    Suppose the detection property \eqref{dp} holds for some collection of subgroups $\mathcal{H}$ and $\cat(H)$ is compactly generated for each $H \in \mathcal{H}$. Denote the inclusion $h: H \to G$ for each $H \in \mathcal{H}$. Then \[\locm(X) = \locm\{\bigoplus\limits_{H \in \mathcal{H}}h_!h^*(X)\}\] and so the natural map $\loc(G) \to \lim\limits_{H \in \mathcal{A}_{\mathcal{H}}(G)}\mathbb{L}(H)$ is an isomorphism.
\end{theorem}
\begin{proof}
    Suppose $\mathcal{L}_G$ is some localising tensor ideal of $\cat(G)$. We need to show that $X \in \mathcal{L}_G$ if and only if $h_!h^*(X) \in \mathcal{L}_G$ for every $H \in \mathcal{H}$.
    \par 
    We want to apply \Cref{Alem} for $A = \bigoplus\limits_{H \in \mathcal{H}} h_!(\mathbbm{1}_H)$ and so first we must show that $\locm(A) = \cat(G)$. To show this we argue that any localising tensor ideal containing $A$ also contains the set $\mathcal{G}$ of compact generators of \Cref{cop2}. It then follows that every localising tensor ideal containing $A$ must indeed by the whole of $\cat(G)$. 
    \par 
    To this end, let $h_!(X)$ be a compact generator as in the statement of \Cref{cop2}, i.e. $X \in \cat(H)$ for some $H$ and it is a compact generator of $\cat(H)$. Since $\locm(A)$ is tensor ideal, we see that $h_!(X) \otimes A \in \locm(A)$. Using that the tensor commutes with coproduct, we see that $h_!(X) \otimes h_!(\mathbbm{1}_H) \in \locm(A)$. 
    \par 
    Now, the projection formula shows that $h_!(X) \otimes h_!(\mathbbm{1}_H) \cong h_!(X \otimes h^*h_!(\mathbbm{1}_H))$ and so $X \otimes h^*h_!(\mathbbm{1}_H) \in h_!^{-1}(\locm(A))$, which is a localising tensor ideal by \Cref{bothdefined}. Using \Cref{summandlem}, we know that $X \cong X \otimes \mathbbm{1}_H$ is a summand of $X \otimes h^*h_!(\mathbbm{1}_H)$ and hence is also in $h_!^{-1}(\locm(A))$. Altogether, this shows that $h_!(X) \in \locm(A)$ as desired. 
    \par 
    Applying \Cref{Alem} then shows that $X \in \mathcal{L}_G$ if and only if $X \otimes (\bigoplus h_!(\mathbbm{1}_H)) \in \mathcal{L}_G$. We have the following isomorphisms. 
    \[X \otimes A = X \otimes (\bigoplus h_!(\mathbbm{1}_H)) \cong \bigoplus (X \otimes h_!(\mathbbm{1}_H)) \cong \bigoplus h_!h^*(X)\]
    The first isomorphism is using that the tensor commutes with coproducts, whilst the second isomorphism is the projection formula. It follows that $X \in \mathcal{L}_G$ if and only if $h_!h^*(X) \in \mathcal{L}_G$ for each $H\in \mathcal{H}$.
    \par 
    The statement about the map being an isomorphism now follows from \Cref{needtobeframe}.
\end{proof}
\subsection{Linking to the spectrum}
We continue to assume that $\cat$ takes values in tensor-triangulated categories and that we have a collection $\mathcal{H}$ of subgroups for which the detection property holds.
\par 
We also assume that the following two conditions hold for every $H \in \mathcal{H}$:
\begin{enumerate}[label=(\roman*)]
    \item $\cat(H)$ is rigidly-compactly generated 
    \item$\cat(H)$ is stratified in the sense of Barthel-Heard-Sanders \cite[Definition~4.4]{BHSstrat}
\end{enumerate}
In particular, this means that the localising tensor ideals of $\cat(H)$ are in one to one correspondence with subsets of the Balmer spectrum of compact objects in $\cat(H)$. Furthermore, this bijection is given by the Balmer-Favi support as defined in \cite[Section~2]{BHSstrat}. We will denote this $\textup{Supp}_H$. For our purposes here, we may take it as a black box.
\par 
We consider the functor $\textup{Spc}: \mathcal{A}_{\mathcal{H}}(G) \to \textup{Set}$ which sends $H$ to the Balmer spectrum of $\cat(H)^d$ and morphisms $f: H \to K$ to $\textup{Spc}(f)$, see \cite{Balmspec}. Note that we are considering the Balmer spectrum simply as a set; we do not consider the topology here. We also have the contravariant powerset functor $\mathcal{P}: \textup{Set} \to \textup{Set}$; this sends morphisms to their inverse image. 
\begin{lemma}\label{functorsareisomorphic}
    There is a natural isomorphism of functors $\mathbb{L} \to \mathcal{P}\circ \textup{Spc}$.
\end{lemma}
\begin{proof}
    Let $f: K \to H$ be a morphism in $\mathcal{A}_{\mathcal{H}}(G)$. We show that the following diagram commutes. 
    \[\begin{tikzcd}[column sep = 3cm]
        \mathbb{L}(H) \arrow{r}{f_!^{-1}} \arrow[swap]{d}{\textup{Supp}_H} & \mathbb{L}(K)\arrow{d}{\textup{Supp}_K} \\
        \mathcal{P}(\textup{Spc}(H)) \arrow[swap]{r}{\textup{Spc}(f^*)^{-1}} & \mathcal{P}(\textup{Spc}(K)) 
    \end{tikzcd}\]
    By the Avrunin-Scott identity \cite[Corollary~4.19]{barthel2023cosupport}, which we may apply since we have assumed $\cat(H)$ and $\cat(K)$ are rigidly-compactly generated, we know that for any $\mathcal{L} \in \loc(Ht)$, we have $\textup{Spc}(f^*)^{-1}(\textup{Supp}_H(\mathcal{L})) = \textup{Supp}_K(f^*(\mathcal{L}))$. Therefore, it suffices to show that $\textup{Supp}_K(f^*(\mathcal{L})) = \textup{Supp}_K(f_!^{-1}(\mathcal{L}))$.
    \par 
    We have shown in \Cref{lreqca} that $\locm(f^*(\mathcal{L})) = f_!^{-1}(\mathcal{L})$, which implies the claim since the Balmer-Favi support classifies the localising ideals in $\cat(K)$ by assumption.
    \par 
    Finally, since $\cat(H)$ and $\cat(K)$ are stratified by assumption, both the vertical maps are isomorphisms. This proves the statement.
\end{proof}
We now reach the main theorem of this section. We will be referring back to it later in the applications and so for convenience we briefly recall all the necessary assumptions. We have a collection of subgroups $\mathcal{H}$ of $G$ such that the detection property holds for $\mathcal{H}$. We also know that $\cat(H)$ is rigidly compactly generated and stratified for each $H \in \mathcal{H}$. Given this, we have the following.
\begin{theorem}\label{biggerboy}
    There is a one to one correspondence between $\loc(G)$ and subsets of $\underset{{H \in \mathcal{A}_{\mathcal{H}}(G)}}{\textup{colim}}\textup{Spc}(\cat(H)^d)$.
\end{theorem}
\begin{proof}
    It follows from \Cref{functorsareisomorphic} that 
    \[\lim\limits_{H \in \mathcal{A}_{\mathcal{H}}(G)}\mathbb{L}(H) \cong \lim\limits_{H \in \mathcal{A}_{\mathcal{H}}(G)}\mathcal{P}(\textup{Spc}(\cat(H)^d))\]
    The powerset functor is contravariant and representable and so it sends colimits to limits; this means that $\lim\limits_{H \in \mathcal{A}_{\mathcal{H}}(G)}\mathbb{L}(H)$ can be identified with subsets of $\underset{{H \in \mathcal{A}_{\mathcal{H}}(G)}}{\textup{colim}}\textup{Spc}(\cat(H)^d)$. Finally, this in turn is isomorphic to the set $\loc(G)$ by \Cref{projloclem}.
\end{proof}
\subsection{An isomorphism of frames}\label{framsection}
We continue to work under the same assumptions as the previous section; in particular, $\cat$ takes values in tensor-triangulated categories and we have a set of localising tensor ideals in $\cat(G)$. In fact, this is actually a poset under inclusion. Our aim in this section is to show that we can determine its structure as a frame.
\par 
In order to do this, we recall the following definitions; see, for example, \cite{KLWell}. A poset $(P, \leq)$ is called a lattice if all non-empty finite joins and non-empty finite meets exist; we will denote the join and meet by $\vee$ and $\wedge$ respectively. Furthermore, we say that $P$ is a frame if all (small) joins exist and the following distributivity law is satisfied for any set $I$ and any $x, \{y_i\}_{i \in I}$:
\begin{equation}
    x \wedge(\bigvee\limits_{i \in I}y_i) \leq \bigvee\limits_{i \in I}(x \wedge y_i)
\end{equation}
A frame homomorphism is then a poset homomorphism which preserves finite meets and arbitrary joins. 
\par 
A localising tensor ideal $\mathcal{L}$ is called radical if $X^{\otimes n} \in \mathcal{L}$ for some $n \in \mathbb{N}$ implies that $X \in \mathcal{L}$.
\par 
As in \cite[Section~5]{KLWell}, we know that $\loc(G)$ is a lattice where the meet of a subset is given by their intersection and the join is the smallest localising tensor ideal containing the union.
\par 
We now want to show that the isomorphism in \Cref{needtobeframe} is actually an isomorphism of frames. For this, it is important to note that by \cite[Lemma~3.14]{barthel2023descenttensortriangulargeometry} we know that the forgetful functor from frames to sets creates all limits, and hence the limit of sets in \Cref{needtobeframe} is also a limit in the category of frames.
\begin{lemma}\label{indisframe}
    Suppose $h:H \to G$ is a faithful morphism. Then $h_!^{-1}$ is a frame morphism.
\end{lemma}
\begin{proof}
    Let $\{\mathcal{L}_i\}_{i \in I}$ be a collection of objects in $\loc(G)$. Their meet is simply given by the intersection, in which case it is clear that $h_!^{-1}(\bigwedge\limits_{i\in I}\mathcal{L}_i) = \bigwedge\limits_{i\in I}h_!^{-1}(\mathcal{L}_i)$.
    \par 
    The more difficult part is showing that $h_!^{-1}$ preserves joins. Explicitly, we need 
    \[h_!^{-1}(\locm(\bigcup\limits_{i \in I}\mathcal{L}_i)) = \locm(\bigcup\limits_{i \in I}h_!^{-1}(\mathcal{L}_i))\]
    Note first that $\bigcup\limits_{i \in I}h_!^{-1}(\mathcal{L}_i) \subseteq h_!^{-1}(\locm(\bigcup\limits_{i \in I}\mathcal{L}_i))$, and since the latter is a localising tensor ideal by \Cref{bothdefined} it follows that $\locm(\bigcup\limits_{i \in I}h_!^{-1}(\mathcal{L}_i)) \subseteq h_!^{-1}(\locm(\bigcup\limits_{i \in I}\mathcal{L}_i))$. 
    \par 
    For the reverse inclusion, suppose $X \in \cat(H)$ is such that $h_!(X) \in \locm(\bigcup\limits_{i \in I}\mathcal{L}_i)$. Now, consider the full subcategory $\mathcal{S}$ of $\cat(G)$ of all objects $Y$ such that $h^*(Y) \in \locm(\bigcup\limits_{i \in I}h_!^{-1}(\mathcal{L}_i))$. This is clearly a localising tensor ideal. Furthermore, for any $Y_i \in \mathcal{L}_i$, we know by the projection formula that $h_!h^*(Y_i) \in \mathcal{L}_i$ and hence $Y_i \in \mathcal{S}$. It follows that $\locm(\bigcup\limits_{i \in I}\mathcal{L}_i) \subseteq \mathcal{S}$. 
    \par 
    Now, we see that $h_!(X) \in \mathcal{S}$ and so $h^*h_!(X) \in \locm(\bigcup\limits_{i \in I}h_!^{-1}(\mathcal{L}_i))$. Since $X$ is a summand of $h^*h_!(X)$ by \Cref{summandlem}, it follows that $X \in \locm(\bigcup\limits_{i \in I}h_!^{-1}(\mathcal{L}_i))$ as we wished.
\end{proof}
\begin{proposition}
    Suppose the detection property \eqref{dp} holds for a collection of subgroups $\mathcal{H}$ and every $\mathcal{L}_H \in \loc(\cat(H))$ is radical for all $H \in \mathcal{H}$. Then every localising tensor ideal of $\cat(G)$ is radical. Consequently, the lattice of localising tensor ideals of $\cat(G)$ is a frame.
\end{proposition}
\begin{proof}
    Suppose $\mathcal{L}$ is a localising tensor ideal and $X\in \cat(G)$ is such that $X \otimes X \in \mathcal{L}$. For every $H \in \mathcal{H}$, applying \Cref{profor} shows that we have the following isomorphism 
    $h_!h^*(X \otimes X) \cong X \otimes X \otimes h_!(\mathbbm{1}_H)$. Since $\mathcal{L}$ is tensor ideal, this must also be contained in $\mathcal{L}$.
    \par This shows that $h^*(X)\otimes h^*(X) \in h_!^{-1}(\mathcal{L})$, which we have assumed is radical. Therefore, $h^*(X) \in h_!^{-1}\mathcal{L}$ and so $h_!h^*(X) \in \mathcal{L}$. This holds for each $H \in \mathcal{H}$ and so from \Cref{projloclem} we see that $X \in \mathcal{L}$. A simple induction shows that if $X^{\otimes n} \in \mathcal{L}$ for any $n \in \mathbb{N}$ then $X \in \mathcal{L}$. i.e. $\mathcal{L}$ is radical. The final statement now follows from \cite[Lemma~5.2]{KLWell}.
\end{proof}
\begin{theorem}
    Suppose the detection property \eqref{dp} holds for a collection of subgroups $\mathcal{H}$ of $G$. Then the natural map $\loc(G) \to \lim\limits_{H \in \mathcal{A}_{\mathcal{H}}(G)}\mathbb{L}(H)$ is an isomorphism of frames.
\end{theorem}
\begin{proof}
    We need to show that $\bigwedge (h^*)^{-1}$ is a frame morphism; the result then follows from \Cref{needtobeframe} and \Cref{indisframe}. It clearly preserves meets, since these are given by intersections. 
    \par 
    Suppose that $\{(\mathcal{L}^i_H)\}_{i \in I}$ is a collection of objects in $\lim\limits_{H \in \mathcal{A}_{\mathcal{H}}(G)}\mathbb{L}(H)$. We need to show the following equality. 
    \begin{equation}\label{ohno}\bigwedge\limits_{H \in \mathcal{H}}(h^*)^{-1}(\locm(\bigcup\limits_{i \in I}(\mathcal{L}_H^i))) = \locm(\bigcup\limits_{i\in I}(\bigwedge\limits_{H \in \mathcal{H}}(h^*)^{-1}(\mathcal{L}_H^i)))\end{equation}
    Fix some $j \in I$ and suppose we have $X \in \cat(G)$ such that $h^*(X) \in \mathcal{L}_H^j$ for all $H \in \mathcal{H}$. Then clearly $h^*(X) \in \locm(\bigcup\limits_{i \in I}(\mathcal{L}_H^i))$. This exactly means that $X$ is in the left hand side of \eqref{ohno}. This holds for any $X \in \bigcup\limits_{i\in I}(\bigwedge\limits_{H \in \mathcal{H}}(h^*)^{-1}(\mathcal{L}_H^i))$ and hence we have an inclusion 
    \[\bigcup\limits_{i\in I}(\bigwedge\limits_{H \in \mathcal{H}}(h^*)^{-1}(\mathcal{L}_H^i)) \subseteq \bigwedge\limits_{H \in \mathcal{H}}(h^*)^{-1}(\locm(\bigcup\limits_{i \in I}(\mathcal{L}_H^i)))\]
    However, from \Cref{bothdefined} it is clear that the right hand side is in fact a localising tensor ideal, and so we have the $\supseteq$ inclusion of \eqref{ohno}. 
    \par 
    We now turn to the other inclusion in \eqref{ohno}. We choose any $j \in I$ and $K \in \mathcal{H}$. From \Cref{rightinverse} we know that $k_!^{-1}\bigwedge\limits_{H \in \mathcal{H}}(h^*)^{-1}(\mathcal{L}_H^j) = \mathcal{L}_K^j$. In particular, we see that 
    \[\mathcal{L}_K^j \in k^{-1}_!( \locm(\bigcup\limits_{i\in I}(\bigwedge\limits_{H \in \mathcal{H}}(h^*)^{-1}(\mathcal{L}_H^i))))\]
    Using \Cref{bothdefined} we know that this is a localising tensor ideal, and since it is true for every $j \in I$, we have that 
    \begin{equation}\label{uhoh}\locm(\bigcup\limits_{i \in I}(\mathcal{L}_K^i)) \subseteq k_!^{-1}(\locm(\bigcup\limits_{i\in I}(\bigwedge\limits_{H \in \mathcal{H}}(h^*)^{-1}(\mathcal{L}_H^i))))\end{equation}
    Now, we take any $X \in \cat(G)$ such that $h^*(X) \in \locm(\bigcup\limits_{i \in I}(\mathcal{L}_H^i))$ for every $H \in \mathcal{H}$, i.e. $X$ is in the left hand side of \eqref{ohno}. From \eqref{uhoh} we find that $h_!h^*(X) \in \locm(\bigcup\limits_{i\in I}(\bigwedge\limits_{H \in \mathcal{H}}(h^*)^{-1}(\mathcal{L}_H^i)))$ for every $H \in \mathcal{H}$. Now, during the proof of \Cref{projloclem} we showed that this implies that $X \in \locm(\bigcup\limits_{i\in I}(\bigwedge\limits_{H \in \mathcal{H}}(h^*)^{-1}(\mathcal{L}_H^i)))$ and hence we have shown that $\subseteq$ inclusion of \eqref{ohno}, completing the proof.
\end{proof}
\section{Colocalising, hom-closed subcategories}
We now consider the dual question of classifying the colocalising, hom-closed subcategories. For this, we fix a partial Green 2-functor $\cat$, and so in particular we have internal homs in the categories $\cat(G)$. We will denote these by $\textup{hom}_G(-,-)$. 
\begin{definition}
    Suppose $\mathcal{T}$ is a closed symmetric monoidal additive category with products. Let $\textup{hom}(-,-)$ be the internal hom. We say that a full subcategory $\mathcal{C}$ of $\mathcal{T}$ is: 
    \begin{enumerate}[label=(\roman*)]
        \item colocalising if it is thick and closed under products 
        \item hom-closed if for all $C \in \mathcal{C}$ and $X \in \mathcal{T}$, we have $\textup{hom}(X,C) \in \mathcal{C}$
    \end{enumerate}
    If $\mathcal{T}$ is triangulated, we also require $\mathcal{C}$ to be a triangulated subcategory.
\end{definition}
We let $\coloc(\mathcal{T})$ be the collection of colocalising hom-closed subcategories of $\mathcal{T}$. Similarly, for any $X \in \mathcal{T}$ we will write $\colocm(X)$ to be the smallest colocalising hom-closed subcategory containing $X$. For a faithful functor between groupoids $h: H \to G$, we define $h_*^{-1}$ in an analogous way to $(h^*)^{-1}$. 
\begin{proposition}\label{otherbothdefined}
    Suppose $h: H \to G$ is a faithful morphism of groupoids. Then $(h^*)^{-1}:\coloc(H) \to \coloc(G)$ is well defined, as is $h_*^{-1}: \coloc(G) \to \coloc(H)$.
\end{proposition}
\begin{proof}
    Since $h^*$ is a closed functor (i.e. preserves the internal hom) and is a right adjoint,  and so preserves products, this statement is straightforward. We focus on the statement for $h_*$. In fact, since it preserves products, the only tricky part is to show that $h_*^{-1}(\mathcal{C})$ is hom-closed, if $\mathcal{C}$ itself is hom-closed. Therefore, take $X \in h_*^{-1}(\mathcal{C})$ and $Y \in \cat(H)$ to be arbitrary. 
    \par 
    For any object $Z \in \cat(G)$ we have the following natural isomorphisms. 
    \begin{align}
        \cat(G)(Z, \textup{hom}_G(h_!(Y),h_*(X))) &\cong \cat(G)(Z \otimes h_!(Y),h_*(X)) \\ &\cong \cat(H)(h^*(Z) \otimes h^*h_!(Y), X) \\ & \cong \cat(H)(h^*(Z),\textup{hom}_G(h^*h_!(Y),X)) \\ &\cong \cat(G)(Z,h_*(\textup{hom}_G(h^*h_!(Y),X)))
    \end{align}
    By the Yoneda Lemma, we see that $\textup{hom}_G(h_!(Y),h_*(X)) \cong h_*(\textup{hom}_G(h^*h_!(Y),X))$. Since we have assumed $h_*(X) \in \mathcal{C}$, it follows from this that $\textup{hom}_G(h^*h_!(Y),X) \in h_*^{-1}(\mathcal{C})$. 
    \par 
    Finally, from know that $Y$ is a summand of $h^*h_!(Y)$ from \Cref{summandlem} and hence $\textup{hom}_G(Y,X)$ is a summand of $\textup{hom}_G(h^*h_!(Y),X)$ and so must itself be in $\mathcal{C}$ as desired. 
\end{proof}
We will also need the following identity. 
\begin{lemma}\label{cole}
    Let $h:H \to G$ be a faithful morphism and take $X \in \cat(G)$. Then there is an isomorphism $h_*h^*(X) \cong \textup{hom}_G(h_!(\mathbbm{1}_H),X)$.
\end{lemma}
\begin{proof}
    This follows from the Yoneda Lemma and the following natural isomorphisms, for any $Z \in \cat(G)$. 
    \begin{align}
      \cat(G)(Z,h_*h^*(X)) &\cong \cat(G)(h_!h^*(Z),X) \\ &\cong \cat(G)(Z \otimes h_!(\mathbbm{1}_H),X) \\ & \cong \cat(G)(Z,\textup{hom}_G(h_!(\mathbbm{1}_H),X)) 
    \end{align}
    Note that the second isomorphism uses \Cref{profor}.
\end{proof}
The following is proved similarly to \Cref{summandlem}, using the right adjoint of the diagonal map $\Delta$ used during that proof. We omit the details. 
\begin{lemma}\label{coindsumlem}
    Suppose $h: H \to G$ is a faithful morphism and take $X \in \cat(H)$. Then $X$ is a summand of $h^*h_*(X)$.
\end{lemma}
We continue with a 1-full, 2-full 2-subcategory $\mathfrak{H}$ of $\mathbf{Gpd}^f_{/G}$ which is closed under isocomma squares as before. 
\par 
Given this, we define a 2-functor $\mathbb{C}: \mathfrak{H} \to \textbf{Set}$ as follows. On objects we take a pair $(H,h)$ and let $\mathbb{C}((H,h))$ be $\coloc(\cat(H))$. 
\par 
Given a 1-morphism $(H,h) \to (K,k)$, i.e. a pair $(f, \gamma_f)$, we let $\mathbb{C}((f,\gamma_f)) = f_*^{-1}$. This is well-defined by \Cref{otherbothdefined}.
\par 
The fact that 2-morphisms are sent to equalities follows just as the analogous statement for $\mathbb{L}$, see \Cref{ldefinition}. 
Now, we have the following dual to \Cref{rightinverse}. 
\begin{proposition}\label{rightinversecoind}
    The map $\bigwedge\limits_{(H,h) \in \mathfrak{H}}(h^{*})^{-1}$ is a right inverse to the natural map $\coloc(G) \to\lim\limits_{(H,h) \in \mathfrak{H}}\mathbb{C}((H,h)) $. 
    Explicitly, let $(\mathcal{C}_H)$ be an element of $\lim\limits_{(H,h) \in \mathfrak{H}}\mathbb{C}((H,h))$ as in \eqref{limite}. Then for any $(K,k) \in \mathfrak{H}$ we have $k_*^{-1}\bigwedge\limits_{(H,h) \in \mathfrak{H}} (h^{*})^{-1}(\mathcal{C}_H) = \mathcal{C}_K$. 
\end{proposition}
\begin{proof}
    Using \Cref{cole} and \Cref{coindsumlem} the argument follows as in \Cref{rightinverse}. 
\end{proof}
With \Cref{detecteddefinition} in mind, it is natural to introduce the following definition. 
\begin{definition}
    We say that colocalising, hom-closed subcategories are detected on restriction to $\mathfrak{H}$ if for every $X \in \cat(G)$ we have 
    \[\colocm(C) = \colocm \{\prod\limits_{(H,h) \in \mathfrak{H}}h_*h^*(X)\}\]
\end{definition}
The following is then the dual of \Cref{bigboy} and is proved in the same way.
\begin{theorem}\label{cainv}
    Suppose colocalising hom-closed subcategories are detected on restriction to $\mathfrak{H}$. Then the natural map $\coloc(G) \to \lim\limits_{(H,h) \in \mathfrak{H}}\mathbb{C}((H,h))$ is invertible.
\end{theorem}
We now turn to the question of when colocalising hom-closed subcategories are detected on restriction to $\mathfrak{H}$. It turns out to be equivalent to the property of localising tensor ideals being detected on $\mathfrak{H}$, at least when our standing assumptions hold. Before we show this, we recall the following notation for orthogonal subcategories. 
\begin{definition}\label{orthdef}
Let $\mathcal{S}$ be a subcategory of $\cat(G)$ for some groupoid $G$. We define the following two orthogonal subcategories. 
\[\mathcal{S}^{\perp} = \{X \in \cat(G)\ |\ \cat(G)(S,X) = 0 \textup{ for all } S \in \mathcal{S}\}\]
and 
\[{}^{\perp}\mathcal{S} = \{X \in \cat(G)\ |\ \cat(G)(X,S) = 0 \textup{ for all } S \in \mathcal{S}\}\]
\end{definition}
It is not hard to check that if $\mathcal{L}$ is a localising tensor ideal then $\mathcal{L}^{\perp}$ is a colocalising hom-closed subcategory. Similarly, for any colocalising hom-closed subcategory $\mathcal{C}$ the orthogonal ${}^{\perp}\mathcal{C}$ is a localising tensor ideal. 
\begin{proposition}\label{clsadc}
    Localising tensor ideals are detected on restriction to $\mathfrak{H}$ if and only if colocalising hom-closed subcategories are detected on restriction to $\mathfrak{H}$.
\end{proposition}
\begin{proof}
    We begin by assuming that localising tensor ideals are detected on restriction to $\mathfrak{H}$. Let $\mathcal{C} \in \coloc(\cat(G))$. From \Cref{cole} it is clear that if $X \in \mathcal{C}$ then also each $h_*h^*(X) \in \mathcal{C}$. 
    \par 
    On the other hand, suppose $h_*h^*(X) \in \mathcal{C}$ for each $(H,h) \in \mathfrak{H}$. Consider the full subcategory $\mathcal{S}$ of $Y \in \cat(G)$ such that $\textup{hom}_G(Y,X) \in \mathcal{C}$. This is clearly a localising subcategory, and it is in fact also tensor ideal. To see this, let $Y \in \mathcal{S}$ and $Z \in \cat(G)$. The isomorphism $\textup{hom}_G(Z \otimes Y,X) \cong \textup{hom}_G(Z,\textup{hom}_G(Y,X))$ along with the fact the $\mathcal{C}$ is hom-closed combine to show that $Z \otimes Y \in \mathcal{S}$. 
    \par 
    By assumption and the isomorphism $h_*h^*(X) \cong \textup{hom}_G(h_!(\mathbbm{1}_H),X)$ of \Cref{coindsumlem}, we know that $\mathcal{S}$ contains $h_!(\mathbbm{1}_H)$ for all $(H,h) \in \mathfrak{H}$. Since localising tensor ideals are detected on restriction to $\mathfrak{H}$, we conclude that $\mathbbm{1}_G \in \mathcal{S}$. However, now we see that $X \cong \textup{hom}_G(\mathbbm{1}_G,X) \in \mathcal{C}$, i.e. we have shown that colocalising hom-closed subcategories are detected on restriction to $\mathfrak{H}$.
    \par 
    We now show the reverse implication, and so assume that colocalising hom-closed subcategories are detected on restriction to $\mathfrak{H}$. Take $\mathcal{L} \in \loc(\cat(G))$. Using \Cref{profor} we see immediately that if $X \in \mathcal{L}$ for some $X \in \cat(G)$, then $h_!h^*(X) \in \mathcal{L}$ for all $(H,h) \in \mathfrak{H}$. 
    \par 
    We therefore need to show that if for some $X \in \cat(G)$ we have $h_!h^*(X) \in \mathcal{L}$ for all $(H,h) \in \mathfrak{H}$, then $X \in \mathcal{L}$. We show this first for $X = \mathbbm{1}_G$, and so assume that $h_!h^*(\mathbbm{1}_G) \in \mathcal{L}$. We claim that $\mathcal{L} = \cat(G)$, which implies that $\mathbbm{1}_G \in \mathcal{L}$. 
    \par 
    Let $\mathcal{J}$ be the localising ideal generated by all $h_!(\mathbbm{1}_H)$; we claim that $\mathcal{J}^{\perp} = 0$. 
    \par 
    Indeed, suppose that $X \in \mathcal{J}^{\perp}$ and so $\cat(G)(h_!(\mathbbm{1}_H),X) = 0$. We have the following isomorphisms.
     \[\cat(G)(h_!(\mathbbm{1}_H),X) \cong \cat(H)(\mathbbm{1}_H,h^*(X)) \cong \cat(H)(h^*(\mathbbm{1}_G),h^*(X)) \cong \cat(G)(\mathbbm{1}_G,h_*h^*(X))\]
     This implies that $h_*h^*(X) \in \locm(\mathbbm{1}_G)^{\perp}$, which is a colocalising hom-closed subcategory. Since colocalising hom-closed subcategories are detected on restriction to $\mathfrak{H}$ we conclude that $X \in \locm(\mathbbm{1}_G)^{\perp}$. 
     \par 
     However, $\locm(\mathbbm{1}_G) = \cat(G)$ and hence $X \cong 0$, showing that $\mathcal{J}^{\perp} = 0$. This implies that $\mathcal{J} = \cat(G)$, i.e. $\mathbbm{1}_G \in \mathcal{J}$. 
     \par 
     We now move to the general case, and so suppose that $\mathcal{L} \in \loc(G)$ and for some $Y \in \cat(G)$ we have $h_!h^*(Y) \in \mathcal{L}$ for all $(H,h) \in \mathfrak{H}$. Consider the full subcategory $\mathcal{S}$ of $Z \in \cat(G)$ such that $Y \otimes Z \in \mathcal{L}$. This is easily seen to be a localising tensor ideal subcategory. From \Cref{profor} we know that $h_!h^*(Y) \cong Y \otimes h_!(\mathbbm{1}_H)$ and so, by assumption, we see that $h_!(\mathbbm{1}_H) \in \mathcal{S}$ for all $(H,h) \in \mathfrak{H}$. We have just shown that necessarily we must have $\mathbbm{1}_G \in \mathcal{S}$ and so $Y \in \mathcal{L}$ as desired.  
\end{proof}
In particular, this means in the situations where we have the detection property for a collection of subgroups $\mathcal{H}$, we get a classification of colocalising hom-closed subcategories as well. We record this in the following. 
\begin{theorem}\label{cascdv}
    Suppose we have a collection of subgroups $\mathcal{H}$ such that, for any $X \in \cat(G)$ we have
    \[\colocm(X) = \colocm\{\prod\limits_{H \in \mathcal{H}}h_*h^*(X)\}\]
    Then the natural map $\coloc(G) \to \lim\limits_{H \in \mathcal{A}_{\mathcal{H}}(G)}\mathbb{C}(H)$ is an isomorphism. This is true, for example, if the detection property holds for $\mathcal{H}$.
\end{theorem}
\begin{proof}
    For such a collection of subgroups $\mathcal{H}$, we can define the subcategory $\mathfrak{H}_{\mathcal{H}}$ as in \Cref{frakde}. Just as we proved in \Cref{NCE}, it follows that colocalising hom-closed subcategories are detected on restriction to $\mathfrak{H}_{\mathcal{H}}$. In particular, this implies that the natural map $\coloc(G) \to \lim\limits_{(H,h) \in \mathfrak{H}_{\mathcal{H}}}\mathbb{C}((H,h))$ is invertible. 
    \par 
    To rewrite the limit over the category $\mathcal{A}_{\mathcal{H}}(G)$, we proceed as before in \Cref{needtobeframe}; essentially all we need to check is that the colocalising hom-closed subcategories of a product of categories can be identified with the product of the colocalising hom-closed subcategories. However, this follows similarly to \Cref{propres}. 
    \par 
    The fact that this is true if the detection property holds for $\mathcal{H}$ follows from \Cref{projloclem} and \Cref{clsadc}.
\end{proof}
We wish to compare the functors $\mathbb{C}$ and $\mathbb{L}$; essentially this is comparing the localising tensor ideals and colocalising hom-closed subcategories of $\cat(G)$. 
\par 
It is shown in \cite[Theorem~7.19, Corollary~7.20]{barthel2023cosupport} that if a rigidly-compactly generated tensor-triangulated category is costratified, then it is also stratified and furthermore there is a bijection between the localising tensor ideals and colocalising hom-closed subcategories. This bijection is realised by sending $\mathcal{L} \mapsto \mathcal{L}^{\perp}$, with inverse being given by $\mathcal{C} \mapsto {}^{\perp}\mathcal{C}$. 

\begin{proposition}\label{yep}
    Assume that $\cat(H)$ is costratified for each $H \in \mathcal{H}$. Then there is a natural isomorphism of functors $\mathbb{C} \to \mathbb{L}$.
\end{proposition}
\begin{proof}
    Suppose $H,K \in \mathcal{H}$ and $f: K \to H$ is a morphism.
    Consider the following diagram.
    \[\begin{tikzcd}[column sep = 3cm]
        \mathbb{C}(H) \arrow{r}{\mathcal{C} \mapsto {}^{\perp}\mathcal{C}}\arrow{d}[swap]{f_*^{-1}} & \mathbb{L}(H) \arrow{d}{f_!^{-1}}\\ \mathbb{C}(K) \arrow{r}[swap]{\mathcal{C} \mapsto {}^{\perp}\mathcal{C}} & \mathbb{L}(K)
    \end{tikzcd}\]
    Our claim is that this commutes, and since the two horizontal arrows are isomorphisms by assumption, this implies that it is a natural isomorphism. Take some $\mathcal{C} \in \coloc(\cat(H))$. We need to show that $f_!^{-1}({}^{\perp}\mathcal{C}) = {}^{\perp}(f_*^{-1}(\mathcal{C}))$. 
    \par 
    Take $X \in f_!^{-1}({}^{\perp}\mathcal{C})$ and $Y \in f_*^{-1}(\mathcal{C})$. Then $0 = \cat(H)(f_!(X),f_*(Y)) \cong \cat(K)(X,f^*f_*(Y))$. We know from \Cref{coindsumlem} that $Y$ is a summand of $f^*f_*(Y)$ and hence $\cat(K)(X,Y)$ is a summand of $\cat(K)(X,f^*f_*(Y)) = 0$. This shows that $X \in {}^{\perp}(f_*^{-1}(\mathcal{C}))$. 
    \par 
    Now, we take some $X' \in {}^{\perp}(f_*^{-1}(\mathcal{C}))$ and $Y' \in \mathcal{C}$. We have $\cat(H)(f_!(X'),Y') \cong \cat(K)(X',f^*(Y'))$. Since $Y' \in \mathcal{C}$, we know from \Cref{cole} that $f_*f^*(Y') \in \mathcal{C}$. In particular, $f^*(Y') \in f_*^{-1}(\mathcal{C})$. By assumption on $X'$ we then see that $\cat(K)(X',f^*(Y')) = 0$. Altogether, this means that $X' \in f_!^{-1}({}^{\perp}\mathcal{C})$ as desired. 
\end{proof}
For the following, we again remind ourselves of our assumptions. We asssume that $\cat(G)$ is a tensor-triangulated category and that we have a collection of subgroups $\mathcal{H}$ of $G$ for which the detection property holds. Furthermore, assume $\cat(H)$ is a rigidly-compactly generated costratified tensor-triangulated category for each $H \in \mathcal{H}$. 
\begin{corollary}\label{biggercor}
 There is a one to one correspondence between the localising tensor ideals and colocalising hom-closed subcategories of $\cat(G)$.
\end{corollary}
\begin{proof}
    We have shown in \Cref{yep} that $\mathbb{C}$ and $\mathbb{L}$ are naturally isomorphic functors and so must have isomorphic limits. However, the limit over $\mathbb{C}$ gives the set of colocalising hom-closed subcategories by \Cref{cascdv}, whereas the limit over $\mathbb{L}$ can be identified with the set of localising tensor ideals by \Cref{needtobeframe}.
\end{proof}
Unwrapping the maps in the above proof shows that the bijection is given by 
\[\mathcal{L} \mapsto \bigwedge\limits_{H \in \mathcal{H}}(h^*)^{-1}((h_!^{-1}(\mathcal{L}))^{\perp}) \]
and
\[\mathcal{C} \mapsto \bigwedge\limits_{H \in \mathcal{H}}(h^*)^{-1}({}^{\perp}(h_*^{-1}(\mathcal{C})))\]
where $\mathcal{L} \in \loc(G)$ and $\mathcal{C} \in \coloc(G)$. 
\par 
However, it turns out that these maps have a simpler description: they are simply the orthogonal maps.
\begin{proposition}\label{permaps}
    Let $\mathcal{L} \in \loc(G)$ and $\mathcal{C} \in \coloc(G)$. Then 
    \[\bigwedge\limits_{H \in \mathcal{H}}(h^*)^{-1}((h_!^{-1}(\mathcal{L}))^{\perp}) = \mathcal{L}^{\perp}\]
    and
    \[\bigwedge\limits_{H \in \mathcal{H}}(h^*)^{-1}({}^{\perp}(h_*^{-1}(\mathcal{C}))) = {}^{\perp}\mathcal{C}\]
\end{proposition}
\begin{proof}
    During the proof of \Cref{yep} we showed that $h_!^{-1}({}^{\perp}\mathcal{C}) = {}^{\perp}(h_*^{-1}(\mathcal{C}))$. In a similar manner, we have that $(h_!^{-1}(\mathcal{L}))^{\perp} = h_*^{-1}(\mathcal{L}^{\perp})$. 
    \par 
    We therefore have $\bigwedge\limits_{H \in \mathcal{H}}(h^*)^{-1}((h_!^{-1}(\mathcal{L}))^{\perp}) = \bigwedge\limits_{H \in \mathcal{H}}(h^*)^{-1}(h_*^{-1}(\mathcal{L}^{\perp}))$. Since $\mathcal{L}^{\perp}$ is a colocalising hom-closed subcategory, it follows from \Cref{cainv} that $\bigwedge\limits_{H \in \mathcal{H}}(h^*)^{-1}(h_*^{-1}(\mathcal{L}^{\perp})) = \mathcal{L}^{\perp}$ as desired. 
    \par 
    The second equality is proved similarly.
\end{proof}
\section{Examples}\label{es}
We now turn to classifying the subcategories in various contexts. Our motivating example is that of the stable module category, although as we show there are other examples where the theory works. 
\par 
We begin by recalling the following types of groups, which we will be interested in.
\begin{definition}\cite{krop} 
    Let $\mathfrak{X}$ be a class of groups and let $\textup{H}_0\mathfrak{X} = \mathfrak{X}$. For every ordinal $\alpha > 0$, we define $\textup{H}_{\alpha}\mathfrak{X}$ inductively. If $\alpha$ is a successor ordinal, we let $\textup{H}_{\alpha}\mathfrak{X}$ be the class of groups which admit a cellular action on a finite dimensional contractible cell complex such that each stabiliser is in $\textup{H}_{\beta}\mathfrak{X}$ for some $\beta < \alpha$. For limit ordinals $\alpha$ let $H_{\alpha}\mathfrak{X} = \bigcup\limits_{\beta<\alpha}H_{\beta}\mathfrak{X}$. We define $\textup{H}\mathfrak{X}$ to be the union of $H_{\alpha}\mathfrak{X}$ over all ordinals $\alpha$. Finally, let $\textup{LH}\mathfrak{X}$ be the class of groups such that every finite subset of elements is contained in an $\textup{H}\mathfrak{X}$ subgroup.
\end{definition}
\begin{definition}\cite[Definition~2.1]{MSstabcat} 
    We say that a group $G$ is of type $\Phi_k$ if any $kG$-module $M$ has finite projective dimension if and only if $M{\downarrow}_F^G$ has finite projective dimension for every finite subgroup $F \leq G$.
\end{definition}
We will write $\textup{LH}\Phi_k$ to denote $\textup{LH}\mathfrak{X}$ groups where $\mathfrak{X}$ is the class of groups of type $\Phi_k$. We will also let $\mathfrak{F}$ be the class of all finite groups.
\subsection{Finite groups}
For finite groups, this gives us an alternative proof of \cite[Theorem~9.7]{repstrat} and \cite[Corollary~1.2]{BIKcol} without relying on Quillen Stratification \cite{quill2}. For this, take $\mathbb{G} = \textup{gpd}^f$ to be the category of finite groupoids along with the faithful functors. We let $G$ be finite and $k$ a field of characteristic $p > 0$. Then $G \mapsto \underline{\textup{Mod}}(kG)$ is a Green 2-functor \cite[Example~10.6]{Green2}. Here, we use $\underline{\textup{Mod}}(kG)$ to denote the stable module category constructed by quotienting out by morphisms which factor through a projective module. 
\par 
We let $\mathcal{E}$ be the collection of elementary abelian $p$-subgroups of $G$.
\begin{proposition}
    There is a one to one correspondence between localising tensor ideals of $\underline{\textup{Mod}}(kG)$, colocalising hom-closed subcategory of $\underline{\textup{Mod}}(kG)$, and subsets of $\underset{{E\in \mathcal{A}_{\mathcal{E}}(G)}}{\textup{colim}}\textup{Proj}(H^*(E,k))$.
\end{proposition}
\begin{proof}
    This follows from \Cref{biggerboy} and \Cref{biggercor}, given that the detection property for the collection of elementary abelian $p$-subgroups holds by Chouinard's theorem. Furthermore, there is a homeomorphism $\textup{Spc}(\underline{\textup{Mod}}(kE)^d) \simeq \textup{Proj}(H^*(E,k))$ by work of Benson, Carlson, and Rickard \cite{bcrthick} and the stable module category is stratified and costratified for elementary abelian $p$-subgroups, see \cite[Theorem~8.1]{repstrat} and \cite[Theorem~11.6]{BIKcol}.
\end{proof}
\subsection{Module category}
Throughout this section, we let $k$ be a field of characteristic $p > 0$ and $G$ be an $\textup{H}\mathfrak{F}$ group. 
\par 
Recall from \Cref{examplemod} that we have a partial Green 2-functor whose values on groups is given by $G \mapsto \textup{Mod}(kG)$. 
We consider localising tensor ideals of $\textup{Mod}(kG)$; by this we mean subcategories as in \Cref{allde} which also satisfy the 2-out-of-3 property for short exact sequences. That is, if two of the modules in a short exact sequence are in the subcategory then so is the third. 
Since restriction and induction are exact functors, it is clear that their inverse images preserve localising tensor ideals in this sense, cf. \Cref{bothdefined}. 
\begin{proposition}\label{modulestr}
    There is a one to one correspondence between localising tensor ideals of $\textup{Mod}(kG)$ and $\lim\limits_{F \in \mathcal{A}_{\mathcal{F}}(G)}\mathbb{L}(F)$.
\end{proposition}
\begin{proof}
    We wish to apply \Cref{needtobeframe}. For this, we need to show for any $kG$-module $X$
    \[\locm(X) = \locm\{\bigoplus_{F \in \mathcal{F}(G)}X{\downarrow}_F^G{\uparrow}_F^G\}\]
    Note that the $\supseteq$ inclusion always holds by the projection formula. Suppose we have a localising tensor ideal $\mathcal{L}$ such that $\bigoplus_{F \in \mathcal{F}(G)}X{\downarrow}_F^G{\uparrow}_F^G \in \mathcal{L}$. 
    \par 
    We show that $X{\downarrow}_H^G{\uparrow}_H^G \in \mathcal{L}$ for each $\textup{H}\mathfrak{F}$ subgroup $H$ by induction on the ordinal $\alpha$ such that $H$ is an $\textup{H}_{\alpha}\mathfrak{F}$ group. Taking $H = G$ then shows the claim. The base case is clear since $\textup{H}_{0}\mathfrak{F} = \mathfrak{F}$. 
    \par 
    Now, suppose $H$ is an $\textup{H}_{\alpha}\mathfrak{F}$ subgroup for $\alpha > 0$. Consider the cellular chain complex $0 \to C_n \to \dots \to C_0 \to k \to 0$ coming from the definition of $\textup{H}\mathfrak{F}$ groups, so that each $C_i$ is a direct sum of modules of the form $k{\uparrow}_K^H$ where $K \leq H$ is an $\textup{H}_{\beta}\mathfrak{F}$ subgroup for some $\beta < \alpha$.
    \par 
    We induce this complex to $G$ and tensor with $X$. This gives a complex of the form $0 \to D_n \to \dots \to D_0 \to X \otimes k{\uparrow}_H^G \to 0$ such that each $D_i$ is the direct sum of modules of the form $X \otimes k{\uparrow}_K^G \cong X{\downarrow}_K^G{\uparrow}_K^G$. By the inductive hypothesis, each $D_i \in \mathcal{L}$ and hence $X{\downarrow}_H^G{\uparrow}_H^G \cong X \otimes k{\uparrow}_H^G \in \mathcal{L}$ as desired. 
   
\end{proof}
Note that it is shown in \cite[Theorem~10.4]{repstrat} that for a finite group $F$, there is a bijection between non-zero localising tensor ideals of $\textup{Mod}(kF)$ and subsets of $\textup{Proj}(H^*(F,k))$. We can use this to refine the description of localising tensor ideals in $\textup{Mod}(kG)$ as follows. 
\begin{corollary}
    There is a one to one correspondence between non-zero localising tensor ideals of $\textup{Mod}(kG)$ and subsets of $\underset{{F \in \mathcal{A}_{\mathcal{F}}(G)}}{\textup{colim}}\textup{Proj}(H^*(F,k))$.
\end{corollary}
\begin{proof}
    Suppose $\mathcal{L}$ is a non-zero localising tensor ideal of $\textup{Mod}(kG)$. First, note that for any non-zero $M \in \mathcal{L}$, we must have that $M{\downarrow}_F^G$ is non-zero for each finite subgroup $F \leq G$. If not, we would have that $M$ is isomorphic to zero as a $k$-vector space, and hence is itself isomorphic to zero. 
    \par 
    We define a functor $\mathbb{L}^0$ just as $\mathbb{L}$, but where we only consider non-zero localising tensor ideals of $\textup{Mod}(kF)$ for each finite subgroup $F \leq G$. From the previous paragraph and \Cref{modulestr} we conclude that there is a one to one correpsondence between non-zero localising tensor ideals of $\textup{Mod}(kG)$ and $\lim\limits_{\mathcal{A}_{\mathcal{F}}(G)}\mathbb{L}^0(F)$. 
    \par 
    Now by using \cite[Theorem~10.4]{repstrat} and arguing just as in \Cref{functorsareisomorphic}, we see that there is a natural isomorphism $\mathbb{L}^0(F) \simeq \mathcal{P}\circ\textup{Spc}$. Finally, the argument of \Cref{biggerboy} now shows that the statement holds.
\end{proof}
\subsection{Stable module category}
Let $k$ be a commutative noetherian ring of finite global dimension, and $G$ an $\textup{LH}\Phi_k$ group. We let $\mathbb{G}^f$ to be the subcategory of $\textbf{Gpd}$ consisting of groupoids which can be written as a coproduct of $\textup{LH}\Phi$ groups along with faithful functors. 
\par 
In \cite{kendall2024stablemodulecategorymodel} we constructed a suitable stable module category for a class of infinite groups, including $\textup{LH}\Phi_k$ groups. We will denote the stable category by $\stabcat(kG)$. This may be constructed as the homotopy category of a suitable model structure on $\textup{Mod}(kG)$. We refer to \cite{kendall2024stablemodulecategorymodel} for the details of the construction, since for this work we may take it as a black box. The following sums up the features that we require.  
\begin{theorem}\label{stableproperties}
    Suppose $k$ is a commutative noetherian ring of finite global dimension and $G$ is an $\textup{LH}\Phi_k$ group. Then the following all hold. 
    \begin{enumerate}[label=(\roman*)]
        \item $\stabcat(kG)$ is a compactly generated tensor triangulated category
        \item For any $kG$-module $M$, we have that $M \cong 0$ in $\stabcat(kG)$ if and only if $M{\downarrow}_F^G \cong 0$ in $\stabcat(kF)$ for every finite subgroup $F \leq G$
        \item Restriction has a left and a right adjoint, given by induction and coinduction respectively
        \item The stable category has products and an internal hom
    \end{enumerate}
\end{theorem}
\begin{proof}
    The first part is shown in \cite[Theorem~1.1]{kendall2024stablemodulecategorymodel}.
    \par 
    We now show $(ii)$. It follows from \cite[Theorem~2.30]{kendall2024stablemodulecategorymodel} that $M \cong 0$ in $\stabcat(kG)$ if and only if $M{\downarrow}_H^G \cong 0$ in $\stabcat(kH)$ for every subgroup $H \leq G$ which is of type $\Phi$. Now it suffices to notice that $M{\downarrow}_H^G \cong 0$ for $H$ of type $\Phi$ if and only if $M{\downarrow}_H^G$ has finite projective dimension, see for example \cite[Proposition~2.13]{kendall2024stablemodulecategorymodel}, which, by definition, is true if and only if $M{\downarrow}_F^G$ has finite projective dimension (and hence is stably isomorphic to zero) for every finite subgroup $F \leq H$.  
    \par 
    Given the description of (acyclic) fibrations and (acyclic) cofibrations in the stable category in \cite[Remark~2.15]{kendall2024stablemodulecategorymodel}, it is straightforward to check that induction-restriction form a Quillen adjunction, as does the adjunction restriction-coinduction. Hence, these induce adjunctions on the stable category. In fact, we can furthermore check that these functors all preserve weak equivalences and hence we do not even need to derive them. 
    \par 
    Finally, the existence of products and an internal hom follows abstractly by Brown representability for compactly generated tensor-triangulated categories. However, we can see from \cite[Theorem~4.3.2]{hovey2007model} that the internal hom is simply the derived functor of $\textup{Hom}_k(-,-)$. Also, since the product preserves weak equivalences we see that it induces the product on the stable category without needing to derive. 
\end{proof}
We need to construct a partial Green 2-functor which gives us the stable module category. 
This may be checked directly as in \cite[Proposition~4.3.7]{Mackey2}. Before we do this, we have the following lemma cf. \cite[Lemma~5.2]{MSstabcat}. 
\begin{lemma}\label{ilem}
    Suppose $G$ is an $\textup{LH}\Phi_k$ group and $f: M \to N$ a morphism in the stable category. Then $f$ is a stable isomorphism if and only if $f$ restricts to a stable isomorphism for each finite subgroup $F \leq G$.
\end{lemma}
\begin{proof}
    We complete $f:M \to N$ to a triangle $M \to N \to C(f) \to \Omega^{-1}(M)$. Then $C(f) \cong 0$ if and only if $f$ is an isomorphism. By \Cref{stableproperties} we know that $C(f) \cong 0$ if and only if $C(f){\downarrow}_F^G \cong 0$ for all finite subgroups $F \leq G$, and so the result follows.
\end{proof}
Note that as a consequence of this, we can identify the rigid objects.
\begin{proposition}
    A module $M$ is rigid in $\stabcat(kG)$ if and only if $M{\downarrow}_F^G$ is rigid in $\stabcat(kF)$ for every finite subgroup $F \leq G$.
\end{proposition}
\begin{proof}
    It is well known that symmetric monoidal functors will preserve rigid objects, see e.g. \cite[Remark~6.2]{ttrum}. On the other hand, by definition, $M$ is rigid if and only if the evaluation map $\textup{hom}(M,k) \otimes N \to \textup{hom}(M,N)$ is an isomorphism. By \Cref{ilem} this is an isomorphism if and only if it is an isomorphism on restriction to every finite subgroup. Hence, if $M{\downarrow}_F^G$ is rigid for every finite subgroup $F \leq G$, we see that $M$ is itself rigid.
\end{proof}
We now construct our partial Green $2$-functor for the stable category.
\begin{proposition}
    The assignment $G \mapsto \stabcat(kG)$ extends to a partial Green 2-functor $\mathbb{G}^{op} \to \textup{ADD}$, where $\mathbb{G}$ denotes the subcategory of $\textbf{Gpd}$ consisting of all groupoids which can be written as a coproduct of $\textup{LH}\Phi_k$ groups.
\end{proposition}
\begin{proof}
By additivity, we can reduce this to the case of one object groupoids, i.e. groups. We have shown in \Cref{stableproperties} that restriction to any subgroup has a left and right adjoint. Furthermore, since the tensor product in the stable category is induced by the usual tensor product of modules, we see that the projection formula holds.
\par 
We are left to check the Beck-Chevalley conditions. It is shown in \cite[Proposition~4.3.7]{Mackey2} that a 2-functor satisfies the left Beck-Chevalley condition if and only if a certain natural transformation is invertible. By \Cref{ilem} we may check invertibility on the finite subgroups, where it holds since the stable module category gives a Mackey 2-functor on finite groupoids \cite[Example~4.2.6]{Mackey2}. The right Beck-Chevalley condition holds for the same reason.
\end{proof}
\begin{theorem}\label{ohyea}
     Suppose $k$ is a commutative noetherian ring of finite global dimension and $G$ is an $\textup{LH}\Phi_k$ group. Then there are maps between the following which are bijective.\[\left\{\begin{gathered} \text{Localising tensor ideal}\\ \text{subcategories of $\stabcat(kG)$}
\end{gathered}\;
\right\} \xleftrightarrow{\ \raisebox{-.4ex}[0ex][0ex]{$\scriptstyle{\sim}$}\ }\left\{
\begin{gathered}
  \text{Colocalising hom-closed}\\ \text{subcategories of
    $\stabcat(kG)$}
\end{gathered}
\right\}\; \xleftrightarrow{\ \raisebox{-.4ex}[0ex][0ex]{$\scriptstyle{\sim}$}\ }\left\{
\begin{gathered}
  \text{Subsets of}\\ \underset{{F \in \mathcal{A}_{\mathcal{F}}(G)}}{\textup{colim}}\textup{Proj}(H^*(F,k))
\end{gathered}
\right\}\
\]
\end{theorem}
\begin{proof}
    This is an application of \Cref{biggerboy} and \Cref{biggercor} and so we need to check the assumptions therein. The detection property holds for the collection of finite subgroups by \Cref{stableproperties}. Also, for every finite group $F$,  $\stabcat(kF)$ is a stratified and costratified rigidly-compactly generated tensor-triangulated category by \cite[Corollary~10.4]{barthel2023lattices}, which also shows that there is a homeomorphism $\textup{Spc}(\stabcat(kF)^d) \simeq \textup{Proj}(H^*(F,k))$. The result then follows.
\end{proof}
We know from \Cref{permaps} that the map from localising tensor ideals to colocalising hom-closed subcategories is given by $\mathcal{L} \mapsto \mathcal{L}^{\perp}$, with inverse $\mathcal{C} \mapsto {}^{\perp}\mathcal{C}$. 
\par 
We can also understand the map from localising tensor ideals to subsets of $\underset{{F \in \mathcal{A}_{\mathcal{F}}(G)}}{\textup{colim}}\textup{Proj}(H^*(F,k))$, which gives us a notion of support for $\textup{LH}\Phi$ groups which classifies the localising tensor ideals. 
\par 
In fact, it turns out that this notion of support agrees with Benson's definition \cite[Section~15]{bensoninf}, which we now recall.
\par 
For any finite subgroup $F \leq G$, we let $\rho_F^*: \textup{Proj}(H^*(F,k)) \to \underset{{F \in \mathcal{A}_{\mathcal{F}}(G)}}{\textup{colim}}\textup{Proj}(H^*(F,k))$ be the natural map. We also let $\textup{Supp}_F$ denote the support defined for $\stabcat(kF)$ as in \cite[Section~3]{repstrat}. 
\par 
We define the support of a $kG$ module $M$ to be 
\[\textup{Supp}_G(M) = \bigcup\limits_{F \in \mathcal{F}(G)} \rho_F^*(\textup{Supp}_F(M{\downarrow}_F^G))\]
\begin{proposition}
    The map in \Cref{ohyea} between localising tensor ideal subcategories of $\stabcat(kG)$ and subsets of $\underset{{F \in \mathcal{A}_{\mathcal{F}}(G)}}{\textup{colim}}\textup{Proj}(H^*(F,k))$ is given by $\mathcal{L} \mapsto \bigcup\limits_{M \in \mathcal{L}}\textup{Supp}_G(M)$.
\end{proposition}
\begin{proof}
    This is a matter of unwrapping the maps involved in \Cref{biggerboy}. To begin with, we have the bijective map $\loc(G) \to \lim\limits_{F \in \mathcal{A}_{\mathcal{F}}(G)}\loc(F)$ from \Cref{needtobeframe}, which, by applying \Cref{lreqca}, we can see to be the map sending $\mathcal{L}$ to $(\locm(\mathcal{L}{\downarrow}_F^G))$.
    \par 
    Now, under the isomorphism of functors in \Cref{functorsareisomorphic}, this is mapped to $(\textup{Supp}_F(\locm(\mathcal{L}{\downarrow}_F^G)))$, which is an element of $\lim\limits_{F \in \mathcal{A}_{\mathcal{F}}(G)}\mathcal{P}(\textup{Spc}(\stabcat(kF)^d))$. Finally, in \Cref{biggerboy} we use that the contravariant powerset functor is representable to send this to an element of $\underset{{F \in \mathcal{A}_{\mathcal{F}}(G)}}{\textup{colim}}\textup{Proj}(H^*(F,k))$, which we see to be exactly the map we claim in the statement.
\end{proof}
We may define cosupport in an analogous manner and this will classify the colocalising hom-closed subcategories of the stable category.
\subsection{Homotopy category of projectives}
We consider an $\textup{H}\mathfrak{F}$ group $G$ and a field $k$ of characteristic $p > 0$. The usual tensor product of chain complexes, where we tensor over the field $k$, gives $K(\textup{Proj}(kG))$ the structure of a (potentially) non-unital tensor-triangulated category. 
\par
It is unclear if the category $K(\textup{Proj}(kG))$ will have a tensor unit; for finite groups the tensor unit will be the injective resolution of $k$, see \cite[Proposition~5.3]{BKinjec}, but for infinite groups this will clearly not work. In particular, this means that we do not know if $K(\textup{Proj}(kG))$ is a tensor-triangulated category. \par 
\par 
It still makes sense to define localising tensor ideals as before and we can still apply our results by showing that localising tensor ideals are detected on restriction to finite subgroups. Before we show this, we explain why we have a left partial Green 2-functor. 
\par
We can extend the partial Green 2-functor $G \mapsto \textup{Mod}(kG)$ to the category of chain complexes. Since complexes of projectives are closed under restriction, induction, and tensoring we actually have a left partial Green 2-functor $G \mapsto C(\textup{Proj}(kG))$. The class of contractible complexes of projectives is again closed under restriction, induction, and the tensor product and so we make take the quotient to get a left partial Green 2-functor $G \mapsto K(\textup{Proj}(kG))$. 
\par We need the following lemma. 
\begin{lemma}\label{huh}
    Let $X \in K(\textup{Proj}(kG))$, $H \leq G$ and $Y \in K(\textup{Proj}(kH))$. Suppose that $X{\downarrow}_H^G \in \locm(Y)$. Then $X{\downarrow}_H^G{\uparrow}_H^G \in \locm(Y{\uparrow}_H^G)$. 
\end{lemma}
\begin{proof}
    Consider the full subcategory of $K(\textup{Proj}(kH))$ of chain complexes $Z$ such that $Z{\uparrow}_H^G \in \locm(Y{\uparrow}_H^G)$. This is a localising, tensor ideal subcategory by \Cref{bothdefined}. It clearly contains $Y$ and hence contains $\locm(Y)$. We conclude that $X{\downarrow}_H^G{\uparrow}_H^G \in \locm(Y{\uparrow}_H^G)$ as desired. 
\end{proof}
\begin{proposition}
    Let $G$ be an $\textup{H}\mathfrak{F}$ group and let $X \in K(\textup{Proj}(kG))$. Then 
    \[\locm(X) = \locm(\{X{\downarrow}_F^G{\uparrow}_F^G\ |\ F \leq G \textup{ is a finite subgroup}\})\]
\end{proposition}
\begin{proof}
    We proceed by induction on the ordinal $\alpha$ such that $G$ is an $\textup{H}_{\alpha}\mathfrak{F}$ group. For $\alpha = 0$ the result is clearly true. 
    \par
    Let $\alpha > 0$ and suppose the result holds true for all ordinals $\beta < \alpha$. Take $H \leq G$ to be a $\textup{H}_{\beta}\mathfrak{F}$ subgroup, where $\beta < \alpha$. By \Cref{huh} and induction we infer that $X{\downarrow}_H^G{\uparrow}_H^G \in \locm(\{X{\downarrow}_F^G{\uparrow}_F^G\ |\ F \leq G \textup{ is a finite subgroup}\})$. 
    \par 
    Now consider the cellular chain complex $0 \to C_n \to \dots \to C_0 \to k \to 0$ coming from the definition of $\textup{H}\mathfrak{F}$ groups, so that each $C_i$ is a direct sum of modules of the form $k{\uparrow}_H^G$ where $H \leq G$ is an $\textup{H}_{\beta}\mathfrak{F}$ subgroup for some $\beta < \alpha$. We call this chain complex $C_*$. 
    \par 
    We tensor this with $X$ to get an exact sequence of chain complexes $C_* \otimes X$. Note that this will be exact since we are working over a field $k$, so $X$ is flat over $k$ in each degree and we are tensoring over $k$.
    \par 
    We see that $C_i \otimes X$ is a direct sum of chain complexes of the form $X{\downarrow}_H^G{\uparrow}_H^G$ where $H \leq G$ is a $\textup{H}_{\beta}\mathfrak{F}$ subgroup for $\beta < \alpha$, and hence $C_i \otimes X \in \locm(\{X{\downarrow}_F^G{\uparrow}_F^G\ |\ F \leq G \textup{ is a finite subgroup}\})$. 
    \par 
    The complex $X$ consists of projectives, and so the exact sequence $C_* \otimes X$ in fact is split in each degree; this implies we may split it into short exact sequences which give us exact triangles in the homotopy category. Using that localising tensor ideals are in particular triangulated subcategories, we conclude that $X \in \locm(\{X{\downarrow}_F^G{\uparrow}_F^G\ |\ F \leq G \textup{ is a finite subgroup}\})$. 
    \par 
    On the other hand, we know from \cite[Proposition~5.3]{BKinjec} that for a finite group $F$ the category $K(\textup{Proj}(kF))$ has a tensor unit given by the injective resolution of the trivial module $k$; we will denote this $\mathbbm{1}_F$. We then have an isomorphism $X{\downarrow}_F^G{\uparrow}_F^G \cong X \otimes \mathbbm{1}_F{\uparrow}_F^G$. It follows that $X{\downarrow}_F^G{\uparrow}_F^G \in \locm(X)$.  
\end{proof}
\begin{theorem}
    Suppose $G$ is an $\textup{H}\mathfrak{F}$ group and $k$ a field of characteristic $p > 0$. Then the localising tensor ideals of $K(\textup{Proj}(kG))$ are in one to one correspondence with subsets of $\underset{{F \in \mathcal{A}_{\mathcal{F}}(G)}}{\textup{colim}}\textup{Spec}^h(H^*(F,k))$.
\end{theorem}
\begin{proof}
    From \cite[Theorem~10.1 and Theorem~11.4]{repstrat} we know that for a finite group $F$, the category $K(\textup{Proj}(kF))$ is stratified and there is a homeomorphism $\textup{Spc}(K(\textup{Proj}(kF))^d) \simeq \textup{Spec}^h(H^*(F,k))$. 
    \par 
    We can now apply \Cref{needtobeframe} to conclude that the natural map $\loc(K(\textup{Proj}(kG))) \to \lim\limits_{F \in \mathcal{A}_{\mathcal{F}}(G)}\mathbb{L}(F)$ is an isomorphism. We can now argue just as in \Cref{biggerboy} to conclude that the statement holds. Although we assumed in \Cref{biggerboy} that we are working with a tensor-triangulated category such that the detection property holds, we only need for $K(\textup{Proj}(kF))$ to be a stratified rigidly-compactly generated tensor-triangulated category for each finite subgroup $F \leq G$.
\end{proof}

\subsection{Derived category and variations}
There is also a partial Green 2-functor $G \mapsto D(\textup{Mod}(kG))$ to which we can apply the techniques. However, the result is not very interesting since, if $F$ is a finite group, then $D(\textup{Mod}(kF))$ will only have two localising tensor ideals: the zero ideal and the whole category \cite[Proposition~2.2]{repstrat}.
\par 
For the rest of this section we let $k$ be a field of characteristic $p > 0$ and $G$ be a group of type $\Phi_k$.
\par
Balmer and Gallauer in \cite[Remark~4.21]{MR4468990} define the derived category of permutation modules for finite groups. We extend this definition to infinite groups as follows.
\begin{definition}
Let $\textup{DPerm}(kG)$ be the localising subcategory of $K(\textup{Mod}(kG))$ generated by those permutation modules of the form $k{\uparrow}_H^G$ where $H \leq G$ is a finite subgroup. 
\end{definition}
\begin{remark}
    For (pro)finite groups, Balmer and Gallauer define this category by considering permutation modules induced from all open subgroups. In our setting, it is more natural to consider only the finite subgroups. Indeed, as shown in \cite{BG}, for finite groups this category has different interpretations, in particular it is equivalent to the derived category of cohomological Mackey functors \cite[Section~5]{BG}. In view of how Mackey functors are defined for infinite groups in \cite{MPNMack}, where the finite subgroups play a crucial role, we see why we choose to define $\textup{DPerm}(kG)$ as above. 
\end{remark}
\begin{lemma}\label{ndres}
Let $H \leq G$ be a finite subgroup. Then restriction gives a well defined functor $\textup{DPerm}(kG) \to \textup{DPerm}(kH)$. Similarly, induction gives a well defined functor $\textup{DPerm}(kH) \to \textup{DPerm}(kG)$.
\end{lemma}
\begin{proof}
    Consider the full subcategory of $K(\textup{Mod}(kG))$ of objects $X$ such that $X{\downarrow}_H^G \in \textup{DPerm}(kH)$. This is clearly localising and by the Mackey formula, we can see that it contains $k{\uparrow}_F^G$ for any finite subgroup $F \leq G$. Therefore, it contains $\textup{DPerm}(kG)$. 
    \par 
    The argument for induction is similar.
\end{proof}
Restriction also must preserve coproducts, since in the larger category $K(\textup{Mod}(kG))$ restriction has a right adjoint, and coproducts in $\textup{DPerm}(kG)$ are simply given by coproducts in $K(\textup{Mod}(kG))$. 
\begin{proposition}
    The category $\textup{DPerm}(kG)$ is closed under the tensor product.
\end{proposition}
\begin{proof}
    Consider $X \in \textup{DPerm}(kG)$. A standard localising subcategory argument, see \cite[Lemma~2.2]{BENSON2011953}, shows that $X \otimes \loc(\{k{\uparrow}_H^G\ |\ H \leq G\}) \subseteq \loc(\{X \otimes k{\uparrow}_H^G\ |\ H \leq G\})$.
    \par 
    Now, we have the isomorphism $X \otimes k{\uparrow}_H^G \cong X{\downarrow}_H^G{\uparrow}_H^G$ and it follows from \Cref{ndres} that this is contained in $\textup{DPerm}(kG)$. From this, we find that $X \otimes \textup{DPerm}(kG) \subseteq \textup{DPerm}(kG)$ for any $X \in \textup{DPerm}(kG)$ as desired.
\end{proof}
In particular, this implies that $G \mapsto \textup{DPerm}(kG)$ is a partial Green 2-functor. 
\begin{lemma}\label{dpermres}
    Let $X \in \textup{DPerm}(kG)$. Then $X \cong 0$ if and only if $X{\downarrow}_F^G \cong 0$ for all finite subgroups $F \leq G$.
\end{lemma}
\begin{proof}
    Only the backwards implication needs to be shown and so assume that $X{\downarrow}_F^G \cong 0$ for each finite subgroup $F \leq G$. Therefore, by adjunction we see that $\textup{Hom}_{K(\textup{Mod}(kG))}(k{\uparrow}_F^G,X) = 0$. 
    \par 
    By definition, $X$ is in the localising subcategory generated by all such $k{\uparrow}_F^G$ and so we conclude that $X$ must be isomorphic to zero.
\end{proof}
To show the existence of a tensor unit, we restrict to groups with a finite dimensional model for the classifying space for proper actions $\underline{E}G$; note that by \cite[Corollary~2.6]{MSstabcat} we know that such groups are of type $\Phi$. 
\begin{lemma}
    Suppose $G$ has a finite dimensional model for the classifying space for proper actions $\underline{E}G$. Then the cellular chain complex $C_*(\underline{E}G)$ is the tensor unit in $\textup{DPerm}(kG)$.
\end{lemma}
\begin{proof}
    We write $C = C_*(\underline{E}G)$ and $\tilde{C} = \tilde{C}_*(\underline{E}G)$ as the augmented cellular chain complex. The slight complication is that $\tilde{C}$ may not be an element of $\textup{DPerm}(kG)$. Indeed, it has $k$ in degree $-1$. However, if we take $X \in \textup{DPerm}(kG)$ then we see that $X \otimes \tilde{C}$ and $X \otimes C$ are now both in $\textup{DPerm}(kG)$.
    \par A result of Kropholler and Wall \cite[Corollary~1.5]{kropwall} implies that $\tilde{C}$ splits on restriction to every finite subgroup and in particular $\tilde{C}{\downarrow}_F^G \cong 0$ in $\textup{DPerm}(kF)$. Hence, we see that $X \otimes \tilde{C}$ also splits on restriction to every finite subgroup and from \Cref{dpermres} we conclude that $X \otimes \tilde{C} \cong 0$ in $\textup{DPerm}(kG)$. This implies that the augmentation induces an isomorphism $X \otimes C \to X$ in $\textup{DPerm}(kG)$. 
\end{proof}
For groups with a finite dimensional model for $\underline{E}G$, this gives $\textup{DPerm}(kG)$ a tensor-triangulated structure.
\begin{theorem}
    Suppose $k$ is a field of characteristic $p > 0$ and $G$ is a group with a finite dimensional model for $\underline{E}G$. Then the localising tensor ideals of $\textup{DPerm}(kG)$ are in one to one correspondence with subsets of $\underset{{F\in \mathcal{A}_{\mathcal{F}}(G)}}{\textup{colim}}\textup{Spc}(\textup{DPerm}(kF)^d)$.
\end{theorem}
\begin{proof}
    For a finite group $F$ the localising tensor ideals of $\textup{DPerm}(kF)$ have been classified by Balmer and Gallauer in \cite{balmer2022ttgeometry}; in fact, they show that $\textup{DPerm}(kF)$ is stratified by $\textup{Spc}(\textup{DPerm}(kF)^d)$. We have shown the detection property holds above in \Cref{dpermres}, and so the statement now follows from \Cref{biggerboy}.
\end{proof}
\bibliographystyle{amsplain}
\bibliography{infgrpref} 

\end{document}